\documentclass[12pt]{article}
\usepackage{amssymb}
\usepackage{theorem}

\newtheorem{theorem}{Theorem}[section]
\newtheorem{proposition}[theorem]{Proposition}
\newtheorem{corollary}[theorem]{Corollary}
\newtheorem{lemma}[theorem]{Lemma}
{\theorembodyfont{\rmfamily}
\newtheorem{definition}[theorem]{Definition}
\newtheorem{example}[theorem]{Example}
\newtheorem{remark}[theorem]{Remark}
\newtheorem{conjecture}[theorem]{Conjecture}
}

\include{epsf}

\begin{document}

\def\CP{\mathbb{C}{\rm P}}
\def\tr{{\rm tr}\,}
\def\endproof{{$\Box$}}

\title{Polyhedral K{\"a}hler Manifolds}

\author{Dmitri Panov}

\maketitle

\begin{abstract}

In this article we introduce the notion of
{\it Polyhedral K\"ahler} manifolds, 
even dimensional polyhedral manifolds 
with unitary holonomy. We concentrate on the $4$-dimensional case,
prove that such manifolds are smooth complex surfaces,
and classify the singularities of the metric. 
The singularities form a divisor and the residues of the flat 
connection on the complement of the divisor 
give us a system of cohomological equations.
Parabolic version of Kobayshi-Hitchin correspondence of T.~Mochizuki
permits us to characterize  polyhedral K\"ahler metrics 
of non-negative curvature
on $\mathbb CP^2$ with singularities at complex
line arrangements.

\end{abstract}

\tableofcontents
 
\section{Introduction and results}

First, we recall the notion of a {\it polyhedral metric}, 
a {\it polyhedral manifold}, and give some  basic facts about them.
Consider a 
piecewise linear manifold $M^d$ with a fixed simplicial decomposition.
Let $\Delta_i^d$ be the simplices of highest dimension of this
decomposition. Choose a flat metric on every $\Delta_i^d$ in such a way
that every two simplices that have a common face are glued by an
isometry. This gives a metric on $M^d$, which is called {\it
polyhedral}, and $M^d$ is called  a {\it polyhedral manifold}.

For every point $x$ of a polyhedral manifold $M^d$ we canonically 
associate its {\it tangent cone}, i.e., a cone with polyhedral 
metric such that a neighborhood of its origin is isometric to 
a neighborhood of $x$. At the non-singular points of  $M^d$ the 
tangent cone is the Euclidean space $\mathbb R^{d}$. 
A polyhedral metric has no singularities at  faces of codimension $1$,
but may have singularities at faces of codimension $2$. The 
tangent cone of the points in the interior of such faces is
isometric to the direct product of a $2$-cone  and the flat space 
$\mathbb R^{d-2}$. 
The angle of the $2$-cone is called the {\it conical angle} 
at the face.

The singular locus of a polyhedral metric is naturally stratified.
A point  of $M^d$ is called a {\it metric singularity of
codimension at least $k$} if its tangent cone is not isometric to 
the direct product of $\mathbb R^{d-k+1}$ and a
$(k-1)$-dimensional polyhedral cone. The set of all
metric singularities of codimension at least
$k$ is denoted by $M_s^{d-k}$.

The complement to the singular locus of the metric  
is connected and we can consider the holonomy of the metric
on it. This gives us a representation  
$\pi_1 (M^{d}\setminus M_s^{d-2})\to SO(d)$
(we will consider only orientable manifolds).
For a generic choice of a polyhedral metric this representation 
has an everywhere dense image in $SO(d)$. 

In this work we study even dimensional polyhedral manifolds 
$M^{2n}$ whose holonomy group
is contained in a subgroup of $SO(2n)$ conjugate to $U(n)$.
On the complement to the singularities of the metric 
these manifolds have a complex structure $J$ parallel with respect 
to the flat metric and compatible with it. In addition to the unitarity
of the holonomy we impose one condition.  
For any face $F$ of codimension $2$
consider a simplex $\Delta^{2n}$ that contains $F$ in its border. The 
parallel complex structure $J$ defined in the interior 
of $\Delta^{2n}$ naturally extends to the whole  $\Delta^{2n}$
and we say that $F$ {\it has a holomorphic direction} 
if $F$ is a piece of a holomorphic hyperplane with respect to $J$.

\begin{definition}\label{gooddefinition} 
A polyhedral manifold $M^{2n}$ is called a 
{\it Polyhedral K\"ahler manifold} ($PK$ manifold) 
if the holonomy of its metric belongs 
to a subgroup of $SO(2n)$ conjugate to $U(n)$ and every codimension $2$
face with conical angle $2k\pi$, $k\ge 2$ has a holomorphic direction.
\end{definition}

\begin{remark} Codimension $2$
faces with conical angle different form $2k\pi$ automatically
have complex direction ({\it cf}. section 3) so we don't need to 
impose this condition on them. 
If we don't impose the condition on the faces 
with conical angle $2k\pi$, $k\ge 2$ we obtain
Thurston's $(X,G)$-cone-manifolds modeled on $X=\mathbb C^n$ 
with $G$ the group of unitary isometries of $\mathbb C^n$ [Th].
\end{remark}

\begin{remark} In this work
simplicial decompositions 
are used only to  define the class of $PK$ manifolds and play a 
secondary role. We will mostly think about $PK$ manifolds as spaces 
with a specific metric, and will not distinguish  
manifolds that are isometric but have different 
simplicial decompositions.

\end{remark}

A $2$-dimensional oriented polyhedral surface is automatically
K{\"a}hler (since $SO(2)=U(1)$) and complete classification of such
structures is given in [Tr] (we recall this classification in
Section 2). In the rest of this work we deal mostly with
$4$-dimensional polyhedral K{\"a}hler manifolds. In Section 2
several  elementary examples of such manifolds are given.

A polyhedral metric is called  {\it non-negatively curved} if the
conical angle at every singular face of codimension 2 is smaller
than $2\pi$.
The original motivation for our study of $PK$ metrics is due to
the following remark of Anton Petrunin. {\it The holonomy of a
non-negatively curved polyhedral $\mathbb{C}P^n$ preserves a
symplectic form}  (this is a partial case of a vanishing theorem
proved in [Ch]). This means that a non-negatively curved
$\mathbb CP^{n}$ is $PK$, we discuss this subject in Section 2.

{\bf Known examples.}
An explicit example of a non-negatively curved polyhedral
$\mathbb CP^2$ can be provided by K\"uhnel's $9$ vertices triangulation [BK].
This polyhedral $\mathbb CP^2$ can be obtained as a finite isometric
quotient of a flat complex $2$-torus, and the holonomy
of the metric on $\mathbb CP^2$ is finite. In general for any $n$ there
exist a series of polyhedral metrics on $\mathbb CP^n$  that are
obtained as quotients of complex tori ([KTM]).

Couwenberg, Heckman, and Looijenga  [CML] study 
geometric structures that are more general than $PK$ metrics. 
They obtain  constant holomorphic curvature metrics on $\mathbb CP^n$, 
having as  singular locus complex reflection hyperplane arrangements.
Their approach is different, in particular from the very beginning they 
start with a holomorphic manifold. It should be possible to prove 
that in the case  when curvature is zero their metrics are $PK$.

{\bf Acknowledgments.} This article contains a revised part of my 
PhD written at \'Ecole Polythechnique under supervision 
of Maxim Kontsevich and defended in 2005. 
I am very indebted to Maxim for his insight, inspiration, guidance, 
and help during my PhD, many of his ideas are contained in this work. 
The remark of Anton Petrunin started this research, 
Misha Gromov gave me the first example of a $PK$ metric on $\mathbb CP^2$,
discussion with Dima Zvonkine and Christophe Margerin 
helped to go on and explanation of Takuro Mochizuki of some of his result
permitted me to finish this work. I would also like to thank
Eduard Looijenga, Nikita Markarian, Carlos Simpson, 
Misha Verbitski, and Jean Yves Welschinger. Without discussions, 
comments and encouragement of all these 
people this article would not be written. Finally I would 
like to thank members of the geometry group at Imperial 
for the excellent environment. This work is supported
by EPSRC grant EP/E044859/1.

\subsection{Local properties of $PK$ metrics}

A $PK$ manifold has a natural complex structure defined  outside
the singular locus, it is constant in the local flat coordinates. 
We will prove that for a $4$-dimensional $PK$ manifold $M^4$
this complex structure can be extended to the whole manifold.

\begin{definition}\label{defhol} Let $M^4$ be a $PK$ manifold.
{\it Holomorphic functions} on $M^4$ are defined as 
continuous functions on $M^4$ that are holomorphic on
the complement to the singularities $M^4\setminus M^2_s$.

A {\it holomorphic chart} in $M^4$ is an open subset $U$
with an injective map $\varphi=(f,g): U\to \mathbb C^2$ 
with $f$ and $g$ holomorphic as above 
and such that $\varphi(M^2_s\cap U)$ 
is an analytic subset of  
$\varphi(U)\subset\mathbb C^2$. 
\end{definition}

The following theorem justifies this definition by proving that 
{\it holomorphic functions} and {\it holomorphic charts} on $M^4$
define together a genuine complex structure on $M^4$.

\begin{theorem}\label{pkcomplex}
Every point of a $4$-dimensional $PK$ manifold is contained in 
a holomorphic chart. Holomorphic charts form together 
a holomorphic atlas on $M^4$ and induce on it a structure of
a smooth complex surface. The singular set  
$M^2_s$ is a complex curve for the defined holomorphic structure. 
\end{theorem}

\begin{remark} We don't consider any intermediate smooth structure on 
$M^4$ in order to define a complex structure on it. At the same time, 
it is known that $PL$
manifolds of dimension  up to $6$ have a canonical smooth structure.
\end{remark}

The main step in Theorem \ref{pkcomplex} is the construction
of holomorphic charts for the singularities of the $PK$ metric.
A neighborhood of every singularity embeds isometrically into 
its tangent cone, and it is possible to introduce 
on the tangent cone the structure 
of a single holomorphic chart.
We call these cones {\it polyhedral K\"ahler cones}, and denote 
by $C_K^4$. All $PK$ cones have a natural 
{\it holomorphic Euler vector field} $e$
({\it cf}. Section 3.1) 
that acts by dilatations of the metric, and this 
field is crucial for us.

\begin{theorem}\label{th:linearization}
Let $C_K^4$ be a $4$-dimensional
polyhedral K{\"a}hler cone. There exists a homeomorphism
$\varphi: C_K^4 \to \mathbb{C}^2$ holomorphic outside the singularities 
of the cone and satisfying the following property: The Euler field $e$ written
in coordinates $(z,w)$ of $\mathbb{C}^2$ is given by
$$\frac{1}{\alpha} z\frac{\partial}{\partial z} +
\frac{1}{\beta} w\frac{\partial}{\partial w}$$ 
where $\alpha$ and $\beta$ are positive real numbers. 
The image of the singular locus of $C_K^4$ under the map $\varphi$ 
is given by a union of curves $c_1z^{\alpha}=c_2w^{\beta}$.
\end{theorem}

The singularity is called {\it irrational} if 
$\frac{\alpha}{\beta}\in \mathbb{R}\setminus \mathbb{Q}$. In this
case its tangent cone is isometric to 
the direct product of two $2$-cones $C_1\times C_2$ 
with conical angles $2\pi\alpha$ and
$2\pi\beta$.

The singularity is called {\it rational of type
$(p,q,\alpha)$, $p,q\in \mathbb{N}$} if its Euler field is equal to
$e=\frac{p}{\alpha} z\frac{\partial}{\partial z}+
\frac{q}{\alpha}w\frac{\partial}{\partial w}$ in coordinates $(z,w)$.
Here $p$ and $q$ are relatively prime, $p\le q$ and $\alpha$ is a positive
real number. Sometimes, when the choice of $\alpha$ is not important,
we may also say that the singularity is of type $(p,q)$.
Since $e$ acts by dilatation of the metric, it preserves the singular
locus of the metric. Thus, in the neighborhood of $x$ any irreducible component
of the singular locus is a curve given by one of the equations: $cz^q=w^p$,
$c\ne 0$; $z=0$; or $w=0$. All these curves are flat
with respect to the induced $PK$ metric. Each curve
$cz^q=w^p$ has a conical point at the origin with the {\it same} angle $2\pi\alpha$,
the line $z=0$ has conical angle $\frac{2\pi}{p}\alpha$, and the line
$w=0$ has conical angle $\frac{2\pi}{q}\alpha$.

The next theorem gives a description of the set {\it
$S(p,q,\alpha)$} of equivalence  classes of singularities of
type $(p,q,\alpha)$. There is a slight difference between the cases
$(1=p=q)$, $(1=p<q)$, $(1<p<q)$. The triple $(p,q,\alpha)$ does
not determine the singularity uniquely, 
the singularities of a given type form an
infinite-dimensional space.

\begin{theorem}
\label{Thm:Spq} On a $2$-sphere, consider the set of metrics (up
to isometry) of curvature $4$, with area $\frac{\pi\alpha}{pq}$,
and having an arbitrary number of conical points. Moreover if
$p>1$ and $q>1$ we mark two conical points, while if $p=1$, $q>1$
we mark one conical point. This set of metrics on the sphere is in
natural 1-to-1 correspondence with the set $S(p,q,\alpha)$.
\end{theorem}

Consider a $(p,q,\alpha)$ singularity, let
$(2\pi\beta_1,...,2\pi\beta_n)$ be the conical angles of the $PK$
metric at the singular branches $c_kz^q=w^p,\; (c_k\ne 0)$. Let
$2\pi\beta_z$ be the angle at $z=0$ and $2\pi\beta_w$ the angle
at $w=0$.

\begin{theorem}
\label{Thm:alpha=} The following relation holds:

$$\alpha=\frac{pq}{2}\left(\sum_k (\beta_k-1)+\frac{\beta_z-1}{p}+
\frac{\beta_w-1}{q}\right)+\frac{p+q}{2}$$

\end{theorem}

Theorems ~\ref{pkcomplex} -~\ref{Thm:alpha=}  are proven in
Section 3.

\subsection{Flat connection and topological relations}

By Theorem \ref{pkcomplex} every $4$-dimensional $PK$
manifold is a complex surface
and the singular locus of the $PK$ metric forms a complex curve on it.
Further on will denote the surface by $S$ and the complex curve by $\Gamma$.

The $PK$ metric on $S$ defines a flat meromorphic 
connection on the tangent bundle of $S$ with 
first order poles at $\Gamma$. In Section 4 we study
this connection, especially in the neighborhood of singularities 
of complex codimension $2$. We give a list of conditions 
that imply that a connection on the tangent bundle to a surface is
a connection of a $PK$ metric (Theorem \ref{connectionPK}).
Using the residues of the connection we  write down a system of 
topological relations on the pair $(S,\Gamma)$. This is done in Section 5.

Let us fix some notations.
Irreducible components of $\Gamma$ will be denoted by  $\Gamma_j$.  
For every component $\Gamma_j$ we denote by $2\pi\beta_j$ the 
conical angle at $\Gamma_j$, i.e., the angle of a $2$-cone orthogonal
to any nonsingular point of $\Gamma_j$. The singularities of $\Gamma$
that are not normal crossings are denoted by  $x_i$ and their type
is denoted by $(p_i,q_i,\alpha_i)$.

\begin{definition} For any surface $S$ a collection of divisors 
$\Gamma_j$ with positive weights $\beta_j$ is called a weighted 
arrangements. In the case when there exists a $PK$ metric 
on $S$ with singularities at $\Gamma_j$ of angles $2\pi\beta_j$
we call $(\Gamma_j, \beta_j)$ the {\it weighted arrangement of the 
$PK$ metric} or the {\it $PK$ arrangement}. Sometimes we will mean by
weighted arrangement the whole data
 $(\Gamma_j,\beta_j; x_i,p_i,q_i,\alpha_i)$.
\end{definition}

Define two numbers related to the behavior of $\Gamma_j$ in
the neighborhood of $x_i$. Denote by $\tilde d_{ij}$ the number of
branches (local irreducible components) of $\Gamma_j$ at $x_i$.
Additionally let $d_{ij}$  be the number of branches, except counting
branches $z=0$ and $w=0$  with weights $\frac{1}{p}$
and $\frac{1}{q}$. Denote by $B_{jk}$, $j\ne k$, the number of
intersections of curves $\Gamma_j$ and $\Gamma_k$ that
represent the normal crossing singularity and  define
$B_{jj}$ by
$$B_{jj}=-\Gamma_j\cdot \Gamma_j+\sum_i p_iq_i(d_{ij})^2,$$

where $\Gamma_j\cdot \Gamma_j$ is the self-intersection number of $\Gamma_j$.

\begin{theorem}
\label{Thm:chernclasses} Any weighted arrangement
$(\Gamma_j,\beta_j; x_i,p_i,q_i,\alpha_i)$ of a $PK$ metric 
satisfies the following relations:

\begin{equation}\label{bgderive}
\forall j \quad \sum_k B_{jk}(\beta_k-1)=
-2\chi(\Gamma_j)-K_S\cdot\Gamma_j-\sum_i (d_{ij}(p_i+q_i)-2\tilde d_{ij})
\end{equation}
\begin{equation}
-c_1(TS)=K_S=\sum_j (\beta_j-1)[\Gamma_j]\in H_2(S,\mathbb R)
\end{equation}

Moreover in the case when for every $i$ $p_i=q_i=1$ we have the
following expression for the second
Chern class:
\begin{equation}c_2(TS)=\sum_i(\alpha_i-1)^2+
\sum_{j\ne k}\frac{1}{2}B_{jk}(\beta_j-1)(\beta_k-1)
\end{equation}

\end{theorem}

Here $K_S$ is the canonical class of $S$ and $\chi(\Gamma_j)$  is
the Euler characteristic of the normalization of $\Gamma_j$.

These relations have the following nature.  Relation (1) is a
consequence of the Gauss-Bonnet formula applied to the curve
$\Gamma_j$. For every $\Gamma_j$ the sum of the defects of its conical
points is equal to its Euler characteristic.
Relations (2) and (3)  express the Chern classes of
$TS$ in terms of the residues of the flat connection on $TS$ corresponding
to the $PK$ metric [Oh].

Construction of weighted arrangements satisfying Equations
~\ref{Thm:chernclasses} (1-3) is a problem of independent
interest. It can lead to combinatorial questions of the following
type:

{\it Problem. Classify arrangements of $3n$ lines on
$\mathbb{C}P^2$ such that every line intersect other lines exactly
at $n+1$ points.}

It is easy to see that such line arrangements with
weights $\beta_j=\frac{n-1}{n}$ satisfy Equations
~\ref{Thm:chernclasses} (1-3). This problem appeared previously in [Hir]
and the list of two infinite series and several exceptional arrangements
satisfying these conditions was given  (all these arrangements are 
complex reflection arrangements). Additional  combinatorial 
questions (about {\it simplicial} and limit $PK$ arrangements 
with a {\it cusp}) are formulated in Section 5.

\subsection {Reconstruction of non-negatively curved
$PK$ metrics from weighted arrangements}

One of the  main results of this paper is the following theorem.

\begin{theorem} \label{th:existence}
Consider a weighted arrangements of lines
$(L_j,\beta_j)$ in $\mathbb{C}P^2$ satisfying the following 
conditions

\begin{equation} \label{stabin}
\sum_j(\beta_j-1)=-3,\;\;\;\; 0<\beta_j<1,\;\;\;\; 
\sum_i d_{ij}(\beta_j-1)>-2
\end{equation}

Then the following inequality holds:

\begin{equation}\label{cp2gieseker} 
\sum_k(\alpha_i-1)^2-\sum_j\frac{1}{2}(1-\beta_j)^2B_{jj}-
\frac{3}{2}\le 0
\end{equation}

Moreover, if the equality holds then there 
exists a $PK$ metric on $\mathbb CP^2$ with conical 
angles $2\pi\beta_j$ at $L_j$ (i.e., $(L_j, \beta_j)$ is a $PK$ arrangement.)

\end{theorem}

We prove this theorem in Section 7 after recalling (Section 6)
the technique of parabolic  bundles. To every arrangement that 
satisfies condition of the theorem we associate a parabolic 
structure on the pull-back of the tangent bundle of $\mathbb CP^2$
to the blow up of $\mathbb CP^2$ at the multiple points of 
the arrangement. We prove that constructed parabolic bundle is stable 
and calculate its parabolic Chern characters.
Inequality (\ref{cp2gieseker})
is a consequence of Bogomolov-Gieseker inequality [M2] (see also [Li]).
The existence of a $PK$ metric in the case of equality 
follows form the parabolic version 
of Kobayashi-Hitchin correspondence  from [M2] and additional statements
about logarithmic connection the we prove in Section 4.

\section{Examples of polyhedral K{\"a}hler manifolds}

In this section we recall the classification of
$PK$ structures on complex curves [Tr] and give several
examples of polyhedral K{\"a}hler manifolds of higher dimension.

\subsection{Flat metrics on surfaces}

{\bf Structures on a $2$-dimensional polyhedral cone.} 
A $2$-dimensional polyhedral cone is a very simple object but already 
it supports the majority of geometric structures that 
are essential for this work. Let us describe these structures.  
Consider a $2$-cone $C^2$ with conical angle $2\pi\alpha$.  
Note first that the flat
metric on $C^2\setminus 0$ defines a conformal and hence a 
{\it holomorphic structure} on $C^2\setminus 0$. 
Moreover $C^2\setminus 0$ is by-holomorphic
to $\mathbb C^*$, so we can chose a {\it holomorphic coordinate $z$} on it
(defined up to a multiplicative constant). This coordinate can be 
used to extend the complex structure from $C^2\setminus 0$ to $C^2$.
We have a natural action 
of $\mathbb R^*$ on $C^2$ by homotheties, corresponding
vector field can be complexified and we call it {\it Euler field}.
This field is given in the coordinate $z$ by the formula  
$\frac{z}{\alpha}\frac{\partial}{\partial z}$, the imaginary part of 
the field acts by isometries of the cone.
The metric induces a {\it flat meromorphic  connection}
on the tangent bundle to the cone and it is given by
$\nabla=d+(\alpha-1)\frac{dz}{z}$. The multi-valued flat coordinate on the
cone, i.e., a coordinate in which the connection on $C^2$ is
trivial, is given by $z^{\alpha}$.

Next theorem classifies polyhedral metrics on  surfaces.

\begin{theorem}{\bf (Troyanov [Tr])} Consider a complex curve $\Gamma$
of  genus $g$ with pairwise distinct marked points $x_1,...,x_n$. Let
$\alpha_1,...,\alpha_n$  be real positive numbers such that
$\sum (\alpha_i-1)=2g-2$. Then there is a unique
(up to a real multiplication constant)
complete flat metric on $\Gamma$ with conical points of angles
$2\pi\alpha_i$ at $x_i$ whose conformal structure on
 $\Gamma\setminus\{x_1,...,x_n\}$
is the same as of $\Gamma$ itself.
\end{theorem}

For completeness we give here a proof of the theorem.

\begin{lemma} For every real $\beta_1,...,\beta_n$ such that
$\sum_i \beta_i=0$, there exists a unique meromorphic 1-form $\eta$
on $\Gamma$ with  simple poles with residues $\beta_1,...,\beta_n$
at the points $x_1,...,x_n$, having purely imaginary periods
(i.e., for every closed path $\gamma\in\Gamma$ we have
$\int_{\gamma}\eta\in i\mathbb{R}$).
\end{lemma}

{\bf Proof}.  By Dirichlet's theorem there exists a unique (up to a
constant) real harmonic function $f$ on $\Gamma$, satisfying the
equation $\Delta f=\sum_i \beta_i\delta_{x_i}$. This function has
logarithmic poles at $x_1,...,x_n$. The 1-form $\eta$ is then
given by
$\eta(\overrightarrow{u})=df(\overrightarrow{u})+idf(J\overrightarrow{u})$,
where $J$ defines the complex structure on $T\Gamma$.

\hfill $\square$

\paragraph {Proof of Theorem 2.1.}

{\bf Existence.} Let $\xi$ be a holomorphic differential on $\Gamma$
with simple zeros $y_1,...,y_{2g-2}$.  It defines a flat metric
$\xi\otimes \bar \xi$ on $\Gamma$ with conical points of angle
$4\pi$ at the points $y_i$. Denote by $\nabla$ the corresponding
connection. Consider the 1-form $\eta$ on $\Gamma$ with purely
imaginary periods that has residue $-1$ at any point $y_i$ and
residue $(\alpha_j-1)$ at any point $x_j$.

Let us prove that the connection $\nabla+\eta$ on $\Gamma$ is unitary.
Indeed, the holonomy of $\nabla+\eta$ along a closed path $\gamma$
is given by the formula
$${\rm hol}_{\gamma}(\nabla+\eta)=
{\rm hol}_{\gamma}(\nabla)\exp(-\int_{\gamma}\eta)
=\exp(-\int_{\gamma}\eta).$$

The first equality follows from the definition of holonomy and the
second follows from the fact that the holonomy of $\nabla$ is
trivial. The connection $\nabla+\eta$ defines a unique (up to a real
multiplication constant) flat metric on $\Gamma$. In order to
define it one should fix the metric at any point of $\Gamma$
different from $x_j$ and translate it by means of $\nabla+\eta$ to
other points of $\Gamma$. The metric constructed this way has
singularities exactly at $x_j$, and the conical angles at $x_j$ are
defined by the poles of $\eta$.

{\bf Uniqueness.} Suppose that we have two metrics $g_1$ and $g_2$
satisfying the conditions of the theorem. Then the 1-form
$\nabla_{g_1}-\nabla_{g_2}$ should be holomorphic and it should
have purely imaginary periods (since both  $\nabla_{g_1}$ and
$\nabla_{g_2}$ are unitary), i.e., it is identically zero. Thus
$g_1$ and $g_2$ coincide.

\hfill $\square$

\subsection{Polyhedral K{\"a}hler manifolds of higher dimension}

Recall that a polyhedral manifold is called {\em non-negatively
curved} if the conical angles at all faces of codimension $2$ are at
most $2\pi$.

\begin{proposition}\label{positiveCP} Let $M^{2n}$ be a non-negatively curved
manifold that has a second cohomology class $h\in H^2(M^{2n})$
such that $h^n$ is non zero in $H^{2n}(M^{2n})$.
Then the holonomy of $M^{2n}$ is contained in $U(n)$, i.e.,
such a manifold is $PK$.
\end{proposition}

This proposition is a simple corollary of  results of J.~Cheeger
(see [Ch]), which we will now describe. We don't need these results
in full generality; instead, we give a version sufficient for our considerations.

Let $M^n$ be a polyhedral manifold, and let $M_s^{n-2}$ 
be the subset of all
its metric singularities. Denote by $H_{L_2}^i(M^n)$ the space of $L_2$-harmonic
forms on $M^n\setminus M_s^{n-2}$ that are closed and co-closed.

\bigskip

{\bf Theorem A.} $\mathrm{dim}(H_{L_2}^i(M^n))=b^i(M^n)$.

\bigskip

{\bf Theorem B.} Suppose that the manifold $M^n$ is non-negatively curved.
Then every harmonic form $h$ in $H^i(M^n)$ is parallel, i.e., $\nabla h=0$.

\bigskip

\begin{remark} As it is stated in [Ch] the Theorem {\bf B} 
indicates that non-negatively curved polyhedral manifolds
are analogs of smooth Riemannian manifolds 
with nonnegative curvature operator (rather then smooth manifolds 
with non-negative sectional curvature).
\end{remark}

Let us deduce Proposition \ref{positiveCP}  from theorems {\bf A} and
{\bf B}. We need a simple fact from linear algebra.
\begin{lemma}\label{trivial}
Consider an Euclidean space $V^{2n}$
with a non-degenerate 2-form $w$, $w^n\ne 0$. Denote by $S_w$
the subgroup of $SO(2n)$ that preserves $w$. Then the group $S_w$ is
contained in a subgroup of $SO(2n)$ conjugate to $U(n)$.
\end{lemma}
{\bf Proof.} We can find orthonormal coordinates $(x_i,y_i)$ in $V^{2n}$
such that $w=\sum_i a_i dx_i\wedge dy_i$ ($a_i\ne 0$).  
It is easy to see that every element of $SO(2n)$ 
that preserves $w$ preserves the form
$w'=\sum_i dx_i\wedge dy_i$. The stabilizer of $w'$ in $SO(2n)$ is
exactly $U(n)$.

 \hfill $\square$

{\bf Proof of Proposition \ref{positiveCP}.} Let $g$ be a non-negatively
curved polyhedral metric on $M^{2n}$. By Theorem {\bf A}
there exists a harmonic 2-form $w$ on
$M^{2n}$ such that $\int_M w^n\ne 0$.
By Theorem {\bf B} $w$ parallel in the flat metric. 
It has constant rank outside
of the singularities, and since $\int_M w^n\ne 0$, $w$ should be
non-degenerate. The holonomy of $g$ preserves $w$, thus by Lemma
\ref{trivial} the holonomy is contained in a subgroup of $SO(2n)$
conjugate to $U(n)$, i.e., g is  a polyhedral K{\"a}hler metric.

\hfill $\square$

This proposition indicates that it should be difficult 
to construct an explicit simplicial decomposition on 
$\mathbb CP^n$ that defines a non-negatively curved metric. 
All examples of $PK$ metrics on $\mathbb CP^n$ 
that we know come from algebraic geometry, and produce a 
metric without a chosen simplicial decomposition.

\medskip

{\bf Examples of non-negatively curved polyhedral $\mathbb{C}P^n$.}

{\bf Example 1.} Choose any non-negatively curved polyhedral metric on $\mathbb{C}P^1$.
Consider the $n$-th symmetric power ${\rm Sym^n}(\mathbb{C}P^1)$ of
$\mathbb{C}P^1$ with induced polyhedral metric. We have
$ \rm{Sym^n}(\mathbb{C}P^1)\simeq \mathbb{C}P^n$, and it is clear that
the constructed polyhedral metric on  $\mathbb {C}P^n$ is 
non-negatively curved. This is the first nontrivial example of 
a higher-dimensional 
$PK$ manifold that I learned and it was 
proposed to me by M. Gromov.   

For $n=2$ we obtain a $PK$ metric on $\mathbb{C}P^2$ with singularities 
at a conic and several lines  tangent to it. 
The conical angle at the conic is equal to $\pi$ 
and the sum of defects of the conical angles
at the lines is equal to $-4\pi$. The conic is the image of the
diagonal of $\mathbb{C}P^1\times \mathbb{C}P^1$.

{\bf Example 2.} Let $T^2$ be a $2$-torus with a flat metric.
Consider the $(n+1)$-th power of $T^2$,
$T^{2n+2}=(T^2)^{n+1}$. Let $T^{2n}$
be a subtorus of $T^{2n+2}$ given by the equation
$\sum_i x_i=0, x_1,...,x_{n+1}\in T^2$.
Let $S_{n+1}$ be the permutation group acting on $T^{2n+2}$.

\begin{lemma} The quotient $T^{2n}/S_{n+1}$ is a $\mathbb{C}P^n$
with a non-negatively curved polyhedral metric.
\end{lemma}

{\bf Proof.} Let $E$ be the unique elliptic curve with the same
conformal structure as $T^2$. Let $L_n$ be a complex line bundle
over $E$ with first Chern class $n+1$. Then
$T^{2n}/S_{n+1}$ can be identified with the space of zero divisors of
sections of $L_n$.  

 \hfill $\square$

For $n=2$ we obtain a $PK$ metric on
$\mathbb{C}P^2$ singular along an elliptic curve of
degree $6$ with $9$ cusps. This curve is
projectively dual to a smooth cubic.

\medskip

\begin{remark} $S^4$ and $\mathbb CP^2$ are the 
only orientable $4$-dimensional manifolds 
that admit a non-negatively curved polyhedral 
metric with irreducible holonomy. 
The case of $S^2\times S^2$ was treated in [Or] using the 
theory of Alexandrov spaces, the results of this paper 
can be obtained in a different way using complex geometry.
\end{remark}

{\bf $PK$-metrics via branched covering.}

One can construct polyhedral K{\"a}hler metrics via 
branched coverings. Let $f:S_2\to S_1$ be a branched covering of a smooth
complex surface $S_1$ by a smooth complex surface $S_2$. Suppose
that $S_1$ has a polyhedral K{\"a}hler metric and $f$ is ramified
over a set of flat  curves on $S_1$. Then the pull-back
of the metric on $S_2$ is a polyhedral K{\"a}hler metric.

Consider the map $f:\mathbb{C}P^2\to \mathbb{C}P^2$,
$f(x:y:z)=(x^n:y^n:z^n)$. This map is ramified at the lines
$x=0,y=0,z=0$. The following two examples use this map 
to produce new $PK$ metrics. 

{\bf Example 3.} $7$ lines.

Consider a $PK$ metric on $\mathbb{C}P^2$ with the
singular locus given by  the lines $x=0,y=0,z=0$ and a conic tangent to
these lines ({\it cf}. Example 1). The conical angle at the conic  is
$\pi$ and the conical angles at the lines are equal to
$2\pi\alpha$, $2\pi\beta$, $2\pi\gamma$, $\alpha+\beta+\gamma=1$.
Consider the branched covering $f(x:y:z)=(x^2:y^2:z^2)$. Then the
singular locus of the pullback metric is composed of 7 lines, 4 of
which have conical angle $\pi$ and three of which have angles
$4\pi\alpha$, $4\pi\beta$, $4\pi\gamma$.

{\bf Example 4.} A metric on a symmetric $K3$ surface.

Consider a $PK$ metric on
$\mathbb{C}P^2$ with the singular locus given by 
the lines $x=0$, $y=0$, $z=0$,
$x+y=z$ (all of them having conical angle $\pi$)
and a conic tangent to these lines (Example 1). Consider the
pull-back metric on $\mathbb{C}P^2$ under the map
$f(x:y:z)=(x^6:y^6:z^6)$.
The preimage of the line $x+y=z$ is given by the equation
$x^6+y^6=z^6$. The double cover of  $\mathbb{C}P^2$ ramified
over the curve $x^6+y^6=z^6$ is a $K3$ surface. This construction
gives a polyhedral K{\"a}hler metric on it.

{\bf Example 5.} $PK$ metrics on  algebraic Kummer surfaces.

Recall that a Kummer $K3$ surface is obtained from a 
complex torus $T^2$ by the quotient with respect to the involution 
$I: x\to -x$ and successive blow up of $16$ fixed points. If 
we first blow up  the points on $T^2$ fixed by $I$, we get a 
surface that is a double cover of the Kummer surface.  
So in order to get  a $PK$ metric on a Kummer surface 
it will be sufficient to construct
any $I$-invariant $PK$ metric on $T^2$ 
blown up at $16$ invariant points.   

Let $\Gamma$ be a genus $2$ curve, $\sigma$ its hyperelliptic
involution, and $\mathrm{Jac}_2(\Gamma)$ the Jacobian of
degree $2$ line bundles on $\Gamma$.
Let $g$ be a flat metric with conical points on $\Gamma$,
invariant under $\sigma$ (we suppose that the conformal
structure of $g$ is that of $\Gamma$). The metric $g$ induces a
$PK$ metric $\widetilde g$ on the symmetric
square $\mathrm{Sym}^2(\Gamma)$ of $\Gamma$.

Recall that $\mathrm{Sym}^2(\Gamma)$ is naturally isomorphic 
to the blow up of  $\mathrm{Jac}_2(\Gamma)$ at the point corresponding
to the canonical class of $\Gamma$. Moreover the involution 
$\sigma$ on $\Gamma$ induces the involution $I$ on the blown up of 
$\mathrm{Jac}_2(\Gamma)$. Consider the degree $16$ cover of $\mathrm{Sym}^2(\Gamma)$
corresponding to the subgroup 
$(2\mathbb Z)^4\subset H_1(\mathrm{Sym}^2(\Gamma))$. One can check that 
the involution $I$ lifts to this cover and it fixes $16$ exceptional curves.
Moreover, $I$ fixes the lift of $\widetilde g$.
This finishes the construction.

\section{Singularities of polyhedral
K{\"a}hler manifolds in dimension 4}

Starting from this section we deal only with $4$-dimensional $PK$
manifolds. In the next two subsections we will prove Theorems \ref{pkcomplex}
and \ref{th:linearization}. 
Before doing this let us explain 
why existence of holomorphic charts (Definition \ref{defhol}) 
on a polyhedral K\"ahler manifold  $M^4$  implies immediately 
that $M^4$ is a complex surface. 
 
Indeed, suppose that $M^4$ can be covered by holomorphic charts
$(U_{\alpha},\varphi_{\alpha})$. To prove that 
$(U_{\alpha},\varphi_{\alpha})$ is a holomorphic atlas on $M^4$
we need to show that for every 
$\alpha$ and $\beta$ the gluing map 

$$\varphi_{\alpha}\varphi_{\beta}^{-1}: 
\varphi_{\beta}(U_{\alpha}\cap U_{\beta})\to
\varphi_{\alpha}(U_{\alpha}\cap U_{\beta})$$
is holomorphic. By Definition  \ref{defhol} the map
$\varphi_{\alpha}\varphi_{\beta}^{-1}$ is continuous on 
$\varphi_{\beta}(U_{\alpha}\cap U_{\beta})$ and holomorphic 
on the complement to an analytic subset. So by standard results
it is holomorphic on the whole domain 
$\varphi_{\beta}(U_{\alpha}\cap U_{\beta})$. 

\hfill $\square$

The same argument gives us the following lemma.
\begin{lemma} \label{halfextend}
For $i=0,1$ suppose that every point
of $M^4\setminus M_s^i$ has a holomorphic chart. Then the 
space $M^4\setminus M_s^i$ has a well-defined holomorphic structure. 

\end{lemma}

The proof of Theorem \ref{pkcomplex} will be done in $3$ steps.
First we show that every point in $M_s^2\setminus M_s^1$ is 
contained in a holomorphic chart. Then we prove that 
singularities of pure codimension $3$ don't exit, 
i.e, $M_s^1=M_s^0$. And finally for singularities of 
codimension $4$ the existence of a holomorphic chart  
is claimed by Theorem  \ref{th:linearization}. 

\subsection{Complex structure in codimension $4$ and the Euler field.}

\begin{lemma}\label{extend} 
Every point $x\in M_s^2\setminus M_s^1$ 
is contained in a holomorphic chart. 
In particular the space $M^4\setminus M_s^1$ has a well-defined 
holomorphic structure.
\end{lemma}

{\bf Proof.} It is sufficient to prove this lemma 
for tangent cones of points in  $M_s^2\setminus M_s^1$,
i.e. for $PK$ manifolds that are the direct 
products of a $2$-cone $C^2$ and the  Euclidean plane $\mathbb R^2$.

By Definition \ref{gooddefinition} the complex structure
on $\mathbb R^2\times (C^2\setminus 0)$ is constant with respect
to the flat connection of the metric and invariant
with respect to the holonomy around $(\mathbb R^2,0)$. If the
conical angle $2\pi\beta$ of $C^2$ is not divisible by $2\pi$ 
then this holonomy is
nontrivial, it rotates the tangent planes of the horizontal 
fibers $(*,C^2\setminus 0)$ by the angle $2\pi\{\beta\}$. 
So these fibers are
holomorphic with respect to the complex structure. The fibers
$(\mathbb{R}^2,*)$ are orthogonal to $(*,C^2\setminus 0)$ and so they 
are holomorphic, since  $J$ preserves the metric. Thus the
complex structure on $M^4=\mathbb R^2\times (C^2\setminus 0)$ is
given by the product of the natural  complex structures on $\mathbb
R^2$  and $C^2\setminus 0$. Finally we
note that $C^2\setminus 0$ is bi-holomorphic to $\mathbb C^*$,
so there is a coordinate $z$ on $C^2$ holomorphic on $C^2\setminus 0$
and continuous on $C^2$. The coordinate $z$ together 
with a holomorphic coordinate $w$ on $\mathbb R^2$ define the structure 
of a chart on $\mathbb R^2\times C^2$. The existence of the complex 
structure on  $M^4\setminus M_s^1$ follows now from Lemma \ref{halfextend}.
The case $2\pi\beta=2\pi k$, $k\ge 2$ is similar, the holonomy is trivial
this time but by Definition \ref{gooddefinition} the vertical fiber 
$\mathbb R^2\times 0$ has a holomorphic direction.
 
 \hfill $\square$

\begin{remark} We proved that every $2$-face 
of a $4$ dimensional polyhedral K\"ahler manifold that belongs to 
the singular locus has a {\it holomorphic direction}.
We need to impose the condition on the faces with conical angle 
$2k\pi$, $k\ge 2$, in order to be able to extend the complex structure 
on these faces. Indeed, for a degree $k$ ramified
cover of $\mathbb{C}^2$ with a branching of order $k$ over a totally real
two-dimensional plane, the complex structure on the cover can not be 
extended on the branching locus.
\end{remark}







\begin{definition} Let $M^4$ be a $PK$ manifold and
let $U$ be the universal cover of $M^4\setminus M_s^2$.
The {\it enveloping map} $E$ of $M^4$ is defined as a 
locally isometric map $E: U \to \mathbb C^2$. Equivalently 
this map can be seen as a multi-valued map from $M^4$ to $\mathbb C^2$
that is locally isometric outside of $M_s^2$ and has infinite
ramification at $M_s^2$. The image of $M_s^2$ under the map 
is called {\it branching set} $B(E)$ of $E$, it is composed 
of linear holomorphic faces. Note that $B(E)$ is usually
everywhere dense in $\mathbb C^2$ but
in the case when $B(E)$ is closed the restriction map 
$E: E^{-1}(\mathbb C^2\setminus B(E))\to \mathbb C^2\setminus B(E)$ 
is a covering map.
\end{definition}

\begin{proposition}\label{pr1=0}
Any  $PK$ cone $C_K^4$ that is a product
of $\mathbb R$ with a $3$-cone is isometric
to the product  of $\mathbb C$ with a $2$-cone.
So $4$-dimensional 
$PK$ manifolds can not have singularities of pure codimension $3$.

\end{proposition}

{\bf Proof}. Suppose that $C_K^4$ is isometric to $\mathbb R\times P^3$.
Denote by $v$ the constant
vector field on $C_K^4$ tangent to the vertical lines
$(\mathbb{R},*)$. This field is acting on 
$\mathbb R\times (P^3\setminus 0)$ preserving
the complex structure defined by Lemma \ref{extend}.
Consider the field $J(v)$ obtained from $v$ by the complex
rotation and let $v_{\mathbb C}=v+iJ(v)$ be the 
complexification of $v$. The field $v_{\mathbb C}$ is constant 
in the flat holomorphic coordinates on the complement to the singularities.
Moreover, since the singularities of $C_K^4$ are tangent to $v$
and they are holomorphic on $\mathbb R\times (P^3\setminus 0)$,
$v_{\mathbb C}$ is also tangent to the singularities.

Consider now the enveloping map $E$ of $C_K^4$
and let us show that its branching locus is contained 
in a complex line through $E(0)$ (note that the image of the center $0$ 
of $C_K^4$ is well defined). Indeed, the singularities of $C_K^4$ 
are of the form $\mathbb R\times r_i$, where $r_i$ is a singular ray of $P^3$.
The image of $\mathbb R\times r_i$ under $E$ in $\mathbb C^2$ is a 
complex half-line containing $E(0)$ at its boundary. 
At the same time it is clear that the 
field $v_{\mathbb C}$ descends to a constant filed $E(v_{\mathbb C})$ 
on $\mathbb C^2$, and so all half-lines 
of $B(E)$ are contained in the line $L$ through $E(0)$ tangent to 
$E(v_{\mathbb C})$. 

Since $B(E)\subset L$, the map 
$E: E^{-1}(\mathbb C^2\setminus L)\to \mathbb C^2\setminus L$
is a covering map. But the set $E^{-1}(\mathbb C^2\setminus L)$
is also a cover of the complement in $C_K^4$ to all half-planes 
$\mathbb R\times r$
tangent to $v_{\mathbb C}$ (including all singular half-planes).  
We deduce that the last complement is isometric to 
a product of a punctured $2$-cone with $\mathbb C$ 
and the proposition follows.

\hfill $\square$

Let us sketch an alternative proof of this proposition where 
instead of studying the enveloping map
we work directly with the cone $P^3$. Consider the restriction of the 
field $J(v)$ on $0\times P^3=P^3$. The field $J(v)$ is well defined 
on $P^3\setminus 0$, it is preserved by the holonomy of the metric 
and it is tangent to the singular rays of $P^3$. We will show that 
$P^3$ has at most $2$ singular rays.

Let $S^2$ be the unit sphere centered at the origin of $P^3$.
Consider the following function on $S^2$:
$$f:S^2\to [-1,1]; \; f(x)=\mathrm{cos}(\angle(e_r(x), I(v(x)))$$

We claim that the critical values of this 
function must be equal to  $1$ or 
$-1$. Indeed, if $x$ is a nonsingular point and $e_r(x)$ is not tangent 
to $J(v(x))$, then $f$ has nonzero differential at $x$. If $x$
is singular (i.e., $ x\in r_i$ ), then  $J(v(x))$ is tangent to $r_i$, thus
$f(x)=\pm 1$. It follows from Morse theory that $f$ must have exactly 
two critical points. Thus the number of conical points on $S^2$ 
is at most $2$. A further analysis  shows that $S^2$ is either the unit 
sphere or a  sphere with $2$ conical points admitting an isometric $S^1$
action preserving the points. This proves the proposition.

 \hfill $\square$

\begin{example}\label{symp}
 Let $S^2$ be a unit sphere and $p$ and $q$ be two 
points on it. Consider a ramified degree $n$ cover 
of $S^2$ by a sphere $\tilde S^2$ with ramifications of order 
$n$ at $p$ and $q$. Then the pull-back metric to 
$\tilde S^2$ has $2$ conical points of angles 
$2\pi n$.  Let $P^3$ be a cone over $\tilde S^2$ and 
consider the polyhedral cone $\mathbb R\times P^3$.
The holonomy of the metric on $\mathbb R\times P^3$ is 
trivial for all choices of $p$ and  $q$, 
but the cone admits a $PK$ structure only if $p$ and $q$
are opposite points on $S^2$. Otherwise, the line 
$\mathbb R \times 0$ forms the singular locus of codimension $3$.
\end{example}

\begin{definition} Let $C_K^4$ be a polyhedral K\"ahler cone. The group 
$\mathbb R^*$ is acting on $C_K^4$ by dilatations  (since $C_K^4$ is a cone).
It is clear that this action preserves the holomorphic structure on
$C_K^4\setminus 0$ defined by Proposition \ref{pr1=0}. 
So the vector field $e_r$ generating this 
action can be complexified and the obtained holomorphic vector
field is called {\it Euler vector field} and denoted by $e$.
The imaginary part of the field $e$ is called the spherical component
$e_s=J(e_r)$. It is important that $e_s$  is acting on 
$C_K^4$ by {\it isometries}. 

\end{definition}

\begin{example} Consider the cone $C_K^4$ that is the direct
product of two $2$-cones with conical angles
$2\pi\alpha$ and $2\pi\beta$. Choose holomorphic coordinates $z$
and $w$ on each $2$-cone as in the beginning of
Section 2.1. Then the Euler field is given by
$e=\frac{1}{\alpha}z\frac{\partial}{\partial z}
+\frac{1}{\beta}w\frac{\partial}{\partial w}$. 
Let us decompose  $e$ as above 
in the radial and  spherical components
$e=e_r+ie_s$. If $\frac{\alpha}{\beta}\in
\mathbb{Q}$, then all orbits of the field $e_s$ are closed, and the
field $e_s$ generates an action of $S^1$ on $C_K^4$. For
$\frac{\alpha}{\beta}$ irrational, the closure of a generic  orbit
of $e_s$ is a $2$-torus.
\end{example}






\subsection{Linear coordinates on $PK$ cones}

In this subsection we prove Theorem \ref{th:linearization}. In particular,
we introduce a holomorphic chart  in a 
neighborhood of every singularity of codimension $4$. So, this also
finishes the proof of Theorem  \ref{pkcomplex}.

Let $C_K^4$ be a $4$-dimensional $PK$ cone. We will consider 
the isometric action generated by $e_s$ on $C_K^4$ and will
distinguish two case.

(1) {\it Irrational.} There exists at least one
non-closed orbit.

(2) {\it Rational.} All orbits of the action are closed.
\medskip

\begin{proposition}\label{ircone} 
If at least one of the orbit of the $e_s$ 
action on $C_K^4$ is non-closed then $C_K^4$ is isometric to the 
product of two $2$-cones.

\end{proposition}

Consider the group of isometries of $C_K^4$ preserving its
origin. This is a compact Lie group, 
and the field $e_s$ generates its subgroup isomorphic 
to $\mathbb R^1$ (because at least one orbit of the action on $S^3$
is non-closed). The closure of this subgroup in the group 
of isometries is a compact connected Abelian group, i.e., 
a torus of dimension at
least two. Thus we have a faithful action of $T^2$ on $C_K^4$ by
isometries. 

Let us show that the branching locus of the enveloping map $E$ of $C_K^4$
is contained in the union of two orthogonal lines in $\mathbb C^2$.
Indeed, $T^2$ is acting on $C_K^4$ and this action
induces an action of $\mathbb R^2$ 
on $\mathbb C^2$ equivariant with respect to $E$ and fixing 
the point $E(0)$ in $\mathbb C^2$. This action factors through the standard
action of $T^2$ on $\mathbb C^2$ and it leaves invariant two orthogonal 
lines $L_1$ and $L_2$ through $E(0)$. The branching locus of $E$ is 
a union of lines trough $E(0)$ invariant under $T^2$ action.
Thus the map 
$$E: E^{-1}(\mathbb C^2\setminus L_1\cup L_2)\to 
\mathbb C^2\setminus L_1\cup L_2$$ 
is a covering map. It follows that 
$C_K^4$ is isometric to a product of two $2$-cones.

\hfill $\square$

Theorem \ref{th:linearization} holds for $PK$ cones isometric to 
the direct product of two 
$2$-cones (see the example above). 
So the first case of the theorem is proved.
To treat the second case we will study the action of $e_s$ on
the unit sphere $S^3$ of $C_K^4$
($S^3$ it the set of points lying at distance $1$ from the origin). 
This action is isometric 
and we suppose this time that all orbits are closed.

\begin{lemma} \label{s1s3}
Suppose that all orbits of the action of $e_s$ on $S^3$
are closed. Then there exists $\alpha>0$ such that all orbits except at most
two have period $2\pi\alpha$, and the exceptional orbits have
periods $\frac{2\pi\alpha}{p}$, $\frac{2\pi\alpha}{q}$,
where $p$ and $q$ are co-prime numbers.
Moreover the action is conjugate to 
the action $(z,w)\to (e ^{i\theta p}z, e ^{i\theta q}w)$,
$\theta\in\mathbb R/2\pi\mathbb Z$ on the unit sphere in 
$\mathbb C^2$, $|z|^2+|w|^2=1$.

\end{lemma}

{\bf Proof.} This lemma is standard and follows essentially  
form the fact that $e_s$ is acting on $S^3$ by isometries,
we will just indicate the proof. It is sufficient to show that 
the action induces on $S^3$ the structure of a Seifert fibration, 
in particular all orbits apart from a finite number have period $2\pi\alpha$
and the periods of all exceptional orbits divide $2\pi\alpha$.
Then the lemma follows from the classification of Seifert fibrations
on $S^3$.

Let $o$ be an orbit of $e_s$ on $S^3$, denote by
$2\pi\alpha_1$ its length. Consider the flow on $S^3$ generated by
$e_s$ in time $2\pi\alpha_1$. It is identical on $o$ and induces
a self-map on an invariant slice transversal to $o$. 
This self-map of the slice is an isometry and it
has finite period $n$ (otherwise there exist orbits of $e_s$
that are not closed). It follows that the flow generated by
$e_s$ in time $2n\pi\alpha_1$ induces the identity map
on $S^3$ thus the period of every orbit divides  
$2n\pi\alpha_1$. Since every 
orbit has a neighborhood where all other orbits have period 
$2n\pi\alpha_1$, the number of exceptional 
orbits is finite.

\hfill $\square$

Now, we are ready to give the proof of Theorem \ref{th:linearization}
in the second case. Note first that since the field 
$e_s$ defines an $S^1$ action 
on $C_K^4$, the field $e=e_r+ie_s$  defines a holomorphic $\mathbb C^*$
action on $C_K^4\setminus 0$.

Suppose first that all orbits of the action of $S^1$ on $S^3$ have
the same length, i.e., the pair $(p,q)$ from Lemma \ref{s1s3} 
is $(1,1)$. Then the quotient  space
$(C_K^4\setminus 0)/\mathbb{C}^*$ is a complex curve homeomorphic 
to $S^2$, hence it is $\mathbb{C}P^1$. Thus $C_K^4\setminus 0$
is isomorphic to a holomorphic $\mathbb C^*$ fibration over $\mathbb{C}P^1$.
This fibration can be completed in a unique way to a line bundle 
by adding the zero section. The completed line bundle has first 
Chern class $-1$, (indeed, the associate $S^1$ 
bundle is homeomorphic to the Hopf fibration of $S^3$). We conclude
that $C_K^4\setminus 0$ can be identified with
$\mathbb{C}^2\setminus 0$ and the Euler field has the form
$\frac{1}{\alpha} z\frac{\partial}{\partial z}+
\frac{1}{\alpha}w\frac{\partial}{\partial w}$.

Consider now the case $2\le p<q$. Let us reduce it to the case
$(p,q)=(1,1)$. According to Lemma \ref{s1s3} there are two
orbits of the  action of $S^1$ on $S^3$ of lengths
$2\pi\frac{\alpha}{p}$ and  $2\pi\frac{\alpha}{q}$. Consider the
corresponding orbits $O_p$ and $O_q$ of the action of
$\mathbb{C}^*$ on $C_K^4$. It follows from Lemma
\ref{s1s3} that the triple $(C_K^4,O_p,O_q)$ is homeomorphic
to a triple $(\mathbb{C}^2,\mathbb{C}^1,\mathbb{C}^1)$ composed of a
complex plane and two transversal lines. Thus there
exists a unique ramified covering of the cone $C_K^4$ by another
polyhedral K{\"a}hler cone  $\widetilde{C_K^4}$ of  degree $pq$
that has ramifications of orders $p$ over $O_p$ and $q$
over $O_q$. It is easy to see that constructed  cone
$\widetilde{C_K^4}$ has type $(1,1)$ and that there are
holomorphic  coordinates $(x,y)$ on $\widetilde{C_K^4}$ such that
the Euler field equals $\frac{1}{\alpha}x \frac{\partial}{\partial
x}+ \frac{1}{\alpha}y\frac{\partial}{\partial y}$. The holomorphic
coordinates on $C_K^4$ will then be $z=x^p,w=y^q$ and the Euler
field is  $e=\frac{p}{\alpha} z\frac{\partial}{\partial z}+
\frac{q}{\alpha}w\frac{\partial}{\partial w}$.

\hfill $\square$

\begin{definition}The coordinates $z$ and $w$ constructed above
are called {\it linear coordinates} of a $PK$ cone.  A $PK$
cone is called {\it rational of type
$(p,q,\alpha)$} ($p,q\in \mathbb{N}$) if its Euler field is equal to
$e=\frac{p}{\alpha} z\frac{\partial}{\partial z}+
\frac{q}{\alpha}w\frac{\partial}{\partial w}$ in the
linear coordinates.
The number $\alpha$ is called the {\it conical angle} of the cone.
A cone is called {\it irrational of type $(\alpha_1,\alpha_2)$},
$\frac{\alpha_1}{\alpha_2}\in \mathbb{R}/\mathbb{Q}$,
if its Euler field is equal
to $e=\frac{1}{\alpha_1} z\frac{\partial}{\partial z}+
\frac{1}{\alpha_2}w\frac{\partial}{\partial w}$.
\end{definition}

\begin{remark}  For a rational polyhedral K{\"a}hler cone of type
$(p,q,\alpha)$, all the orbits of the Euler field action are
given by the equations $\frac {z^p}{w^q}=\mathrm{const}$. 
These curves are flat
with respect to the $PK$ metric and all of them
(except the curves  $z=0,w=0$)
have the same conical angle at $0$ equal to $2\pi \alpha$.
\end{remark}

The proof of the following corollary is contained the second 
part of the proof of  Theorem \ref{th:linearization}.

\begin{corollary}\label{coverpq}
For a polyhedral K{\"a}hler  cone $C_1$  of
type $(p,q,\alpha)$ there exists a unique cone  $C_2$ of type
$(1,1,\alpha)$ with a holomorphic map $f:C_2\to C_1$ of degree $pq$
that is a local isometry outside the branching locus.
\end{corollary}

Finally, we describe all $PK$ cones whose singular locus is 
a union of two lines in linear coordinates.

\begin{lemma} \label{2linecone}
Let $C_K^4$ be a $PK$ cone with linear 
coordinates $(z,w)$ and such that  the singular locus is the union 
of the lines $z=0$ and $w=0$. Then either $C_K^4$ is isometric to 
the product of two $2$-cones, or the metric on $C_K^4$ is the
pull-back of a constant metric on $\mathbb C^2$ 
under the map $\mathbb C^2\to \mathbb C^2$, 
$(z,w)\to (z^n,w^m)$. 
\end{lemma}  

{\bf Proof.} We can suppose that $C_K^4$ is a rational cone, irrational
case is treated by Proposition \ref{ircone}. 
Using Corollary \ref{coverpq} we can assume 
that the cone is of type $(1,1)$, i.e., conical angles at the lines 
$z=0$ and $w=0$ are both equal $2\pi\alpha$. Now, consider two cases.

1) $\alpha$ is not integer. Fix a non-singular point $x$ in $C_K^4$
and consider the holonomy of the metric based at $x$. This 
holonomy is generated by two commuting operators $H_z$ and $H_w$ 
corresponding to two pathes around lines $z=0$ and $w=0$ 
(both operators  are non-trivial, since $\alpha$ is non integer).
Then on $C_K^4\setminus\{zw=0\}$ we have two holomorphic 
rank $1$ sub-bundle of $TC_K^4$  invariant under parallel translation
and orthogonal at every point. 

Consider the enveloping map $E: C_K^4\to \mathbb C^2$. Invariant 
sub-bundles are mapped by $E$ to constant orthogonal sub-bundles of 
$T\mathbb C^2$. It is clear that the ramification locus 
is composed of two lines through $E(0)$, tangent to one of the constant 
fields. These lines are orthogonal and so $C_K^4$ is a direct product
of two $2$-cones.

2) $\alpha$ is integer. Then it is clear that the enveloping map from 
$C_K^4$ is in fact not multivalued but  is a finite degree ramified 
covering of $\mathbb C^2$ with ramifications of degree $\alpha$ at the
lines $z=0$, $w=0$. Moreover, the images of both lines are lines in 
$\mathbb C^2$ containing $E(0)$, so we are in the second case described
by the lemma.

\hfill $\square$

\subsection{$4_{\mathbb R}$ $PK$ cones and  $2_{\mathbb R}$ spheres
with conical points (proof of Theorems \ref{Thm:Spq},
\ref{Thm:alpha=})}

We start the proof of Theorem \ref{Thm:Spq} and 
associate  to every $PK$ cone of type
$(1,1,\alpha)$ a metric on a sphere $S^2$, of curvature 4,
having conical singularities.

Denote by $S^3$ the unit sphere around the origin of the cone
and by $S^2$ the quotient of $S^3$ by the action of $e_s$.
Locally, outside the singularities, the action of $e_s$ on $S^3$
is isometric to the action of the Euler field on the standard
(nonsingular) unit sphere. Therefore, locally, outside the
singularities, the  quotient metric on  $S^2$ is isometric to the
quotient of  the standard (nonsingular) sphere by $e_s$. The last
quotient obviously has curvature 4. The singularities of the cone 
correspond to the conical points on $S^2$.

\hfill $\square$

\begin{lemma}
\label{lm:hopf}
Let $\Omega$ be a contractible domain on the
standard sphere $S^2$ of curvature 4 (without conical points).
Then, for any positive $l$, there is a unique metric $g$ of
curvature 1  on $\Omega\times S^1$ with the following properties:
All the fibers of the product are geodesics of length $l$;
there is an action of $S^1$ on $\Omega\times S^1$ by isometries;
the quotient metric on $\Omega$ coincides with the original one.
\end{lemma}

{\bf Proof.} Let $\phi:S^3\to S^2$ be the standard Hopf fibration.
The universal cover of $\phi^{-1}(\Omega)$ is diffeomorphic to
$\Omega \times \mathbb{R}$, and $\mathbb{R}$ acts on it by
parallel translations. The quotient of $\Omega\times \mathbb{R}$
by the subgroup $l\mathbb{Z}$ of $\mathbb{R}$ induces on
$\Omega\times S^1$ the metric we are looking for.

\hfill $\square$

There is a natural connection $\nabla$ on the fibration $\Omega\times
\mathbb{R} \to \Omega$. Its horizontal distribution is given by
the planes orthogonal to the fibers. The following lemma is standard
and we omit the proof.

\begin{lemma}\label{holgamma} The holonomy of the connection $\nabla$ along
a closed curve $\gamma\subset \Omega$ is equal to
the parallel translation by  $2{\rm area}(\gamma)$, where ${\rm
area}(\gamma)$ is the algebraic area bounded by $\gamma$.
\end{lemma}

Now, let $S^2$ be a sphere with a metric of curvature $4$ with
conical points. We will associate to it a $PK$ cone of type
$(1,1)$. 

First, we reconstruct the sphere $S^3$ of curvature $1$
(with singularities) that fibers over $S^2$.  
Cut $S^2$ by geodesic segments with
vertices at all the conical points, in order to obtain a
contractible polygon $P$. This polygon can be immersed
into the standard sphere of curvature 4 by the enveloping map.
Consider the fibration over $P$ from Lemma \ref{lm:hopf} with
length $l=2{\rm area}(P)$. The holonomy of the fibration along the
border of $P$ is trivial (by Lemma \ref{holgamma} the circle $S^1$ makes
one full rotation). This means that the original gluing of $P$, which
gives $S^2$ with conical points, can be lifted to a gluing of $P\times
S^1$. To construct such a gluing, we choose a horizontal section
$s$ of $P\times S^1$ over the boundary of $P$ and identify
$(x,s(x))$ with $(y,s(y))$ whenever $x$ and $y$ are identified by
the gluing of $P$. Since the border circle turns once, we obtain
the sphere $S^3$.

Now, consider the space $\mathbb{R}_+\times S^3$ with the metric
$(dr)^2+r(ds)^2$, where $r\in \mathbb{R}_+$ and $(ds)^2$ is the
metric constructed on $S^3$. This space is a $PK$-cone of type
$(1,1)$.

\hfill $\square$

{\bf The general case.} The  $2$-sphere with conical and marked points 
associated to the $PK$ cone is given by the metric 
quotient of the unit sphere in the $PK$ cone
by the action of the field $e_s$. 
The marked points correspond to the multiple orbits.

Let $S^2$ be a sphere of curvature 4 with conical points two
of which, $x$ and $y$, are marked. Let us construct for every $1<p<q$
the corresponding $(p,q)$-cone. Take first the 
$(1,1)$-cone $C$ associated to the sphere $S^2$ constructed above.
Denote by $l_x$ and $l_y$ the preimages of $x,y\subset C$ under the
projection to $S^2$. Consider the ramified covering  of degree $pq$ over
$C$ with the branchings of orders $p$ over $l_x$ and of $q$ over $l_y$. 
This is the $(p,q)$-cone we are looking for.

\hfill $\square$

{\bf Proof of Theorem \ref{Thm:alpha=}.}

Consider the case of  $PK$ cones of type
$(1,1,\alpha)$.
Let $S^2$ be the quotient sphere
associated to the cone. It has $n$ conical points of angles
$2\pi\beta_1,...,2\pi\beta_n$ and its area is given by 
the Gauss-Bonnet formula:
$$\mathrm{area}(S^2)=\frac{1}{4}(4\pi+\sum_i 2\pi(\beta-1))$$

From Lemma \ref{holgamma}
$2\pi\alpha=2\mathrm{area}(S^2)$. This proves the theorem for
$(1,1,\alpha)$ singularities. Singularities of other types are
treated in a same way using Corollary \ref{coverpq}.

\hfill $\square$

\subsection {$PK$ metrics on singular piecewise linear spaces}

Definition \ref{gooddefinition} can be naturally extended to the
following class of $PL$ manifolds with singularities. We call a
$4$-dimensional topological space with a simplicial decomposition a {\it
$PL$-manifold up to codimension $2$} if every $3$-simplex is
a border of exactly two $4$-simplices. A
compatible choice of flat metric on the $4$-simplices 
of such a space defines a polyhedral metric on it.
Obtained metric has singularities only in codimension $2$ and we can
repeat Definition \ref{gooddefinition} saying that this metric is a
$PK$-metric if its holonomy is contained in $U(2)$ and all singular
$2$-faces of conical angles $2\pi k$ ($k\ge 2$) 
have holomorphic directions. 
A space with such a structure is
called a {\it singular $PK$ manifold.}

Most of the theorems of this section about (nonsingular)
$4_{\mathbb R}$ dimensional $PK$ manifolds can be restated for
singular $4_{\mathbb R}$ dimensional $PK$ manifolds. In fact,
these singular manifolds are complex surfaces with isolated
singularities. We will formulate the result but will skip the proof.

\begin{theorem}\label{singsurf} For a singular $PK$ manifold of
dimension $4_{\mathbb R}$ its complex structure defined outside 
singularities can be extended to the whole singular manifold.
Obtained complex space is a complex surface with isolated 
singularities. In the neighborhood of every isolated singularity
there is a natural holomorphic field $e=e_r+ie_s$ 
such that the (real) field $e_r$ acts by dilatation 
of the metric and the (real) field $e_s$ generates an
action of $S^1$ by isometries. 
\end{theorem}

\section{Polyhedral K{\"a}hler metrics via logarithmic connections}

For every $PK$ surface $S$ the $PK$ 
metric induces on $T(S\setminus sing)$ 
a holomorphic, flat, unitary, torsion free connection. 
This connection extends to a meromorphic 
connection on $TS$ with 
first order poles at the singular locus. 
In this section we will write explicit formulas of 
$PK$ connections  in {\it linear coordinates} $z,w$ 
on $2$-dimensional $PK $cones. We also 
give a condition for a unitary connection on the tangent bundle of 
a complex surface that implies that the corresponding metric on $S$ is $PK$.

\subsection {Definitions and first results}
Let $M$ be a complex manifold and $D$ be a normal
crossing divisor. A meromorphic $1$-form $\omega$ on $M$ is called {\it
logarithmic} with respect to $D$ if it is holomorphic on 
$M\setminus D$, and
in a neighborhood of any point of $D$ it can be
represented as
$$\omega=\sum_{i=1}^k f_i\frac{dz_i}{z_i}+\sum_{i=k+1}^n f_i dz_i, $$
where $f_i$ are holomorphic functions, $z_i$ are local coordinates,
and $D$ is given  locally by the equation

$$D=\cup_{i=1}^k\{z_i=0\}$$
The sheaf of logarithmic 1-forms is denoted by $\Omega^1(\mathrm{log}D)$.

Let $E$ be a holomorphic vector bundle over $M$.
A meromorphic connection  $\nabla$ on $E$ is called {\it logarithmic}
(with respect to $D$) if it can be written in local coordinates as

$$\nabla=d+A,$$
where $A$ is a $\Omega^1(\mathrm{log}D)$-valued section of $\mathrm{Aut}(E)$.

For any irreducible component $D_i$ of $D$ we denote by
$\mathrm{Res}_{D_i}(\nabla)$ the {\it residue} of $\nabla$ with respect to
$D_i$, it is a holomorphic section of $\mathrm{Aut}(E)|_{D_i}$.

The following proposition is standard, the proof  can be found in  
Section $4$ of the article of Malgrange in [Bo].

\begin{proposition}\label{malgrange} 
Let $\nabla$ be a flat logarithmic connection 
on $(M,E^k)$ with poles at a normal crossing divisor $D$. Suppose that
all eigenvalues of $Res_{D}(\nabla)$ are contained in $]-1,0]$.  

For a point $x\in D$ chose local coordinates $z_i$ such that $D_i$ is given
by $\cup_i\{z_i=0\}$. Then there exists a neighborhood $U$ of $x$ and 
holomorphic sections $s_1,...,s_k$ of $E^k$ giving a trivialisation of 
$E^k$ over $U$ and such that in this trivialisation $\nabla$  is given by:

$$\nabla=d+\sum_i B_i\frac{dz_i}{z_i}$$

where $B_i$ are constant matrix-valued functions.
\end{proposition}

\begin{definition}
For a complex manifold $M$ and a meromorphic connection $\nabla$ on 
$TM$ its {\it torsion} is a meromorphic section of 
$\Omega^2(M)\otimes TM$. It is given by  the following formula:
$$T(u,v)=\nabla_u(v)-\nabla_v(u),$$
$u,v\in TM$. A connection with zero torsion is called
{\it torsion-free}.
\end{definition}

Now, we will restrict our attention to connections on
the tangent bundles of surfaces. 
Let $S$ be 
a surface with a weighted arrangement of curves
 $(\Gamma_{j},\beta_{j})$ 
and  let $x_1,...,x_k$ be the points of the 
arrangement of multiplicity at least $3$.

\begin{definition} 
We say that a meromorphic connection $\nabla$ on $TS$ 
is {\it partially adapted} to $(\Gamma_j, \beta_j)$ if
$\nabla$  is logarithmic on  $S\setminus\{x_1,...,x_k\}$,
$\mathrm{Res}_{\Gamma_{j}}\nabla$ has eigenvalues
$(\beta_{j}-1,0)$ at $\Gamma_{j}$, and $T\Gamma_{j}\cong
\mathrm{ker}(\mathrm {Res}_{\Gamma_{j}}\nabla)$.
\end{definition}

\begin{lemma} \label{half1}
On every $PK$ surface $S$ the $PK$ connection is partially 
adapted to the weighted arrangement of $S$.

\end{lemma}

{\bf Proof.} The statement of the lemma clearly holds 
at smooth points of the singular locus of $S$, because they 
can be embedded isometrically
in the product of $\mathbb C$ with a $2$-cone. 
At the same time the connection on the
$2$-cone is logarithmic and has residue $\beta-1$ 
where $2\pi \beta$ is the cone angle. 

From the description of double points of the singular locus 
(Lemma \ref{2linecone}) it follows that the connection is also 
logarithmic at these points. Indeed,
the connection on the direct product of two $2$-cones is logarithmic, 
and for the second type of cones that are branched covers 
of $\mathbb C^2$ one can change the metric without changing 
the connection to make these cones also direct products.

\hfill $\square$

\begin{lemma}\label{holtorsion} Let $S$ be a complex surface 
with a weighted arrangement $(\Gamma_j,\beta_j)$, $\beta_j\ne 1$. 
Suppose that $\nabla$ is partially adapted, then its torsion is
holomorphic.

\end{lemma}

{\bf Proof}. By Hartogs theorem it is sufficient to 
show that the torsion of $\nabla$ is
holomorphic outside of the multiple points of the arrangement
$\Gamma_{j}$. So it is sufficient to consider the case of a
connection on $\mathbb{C}^2$ with a pole at the line $z=0$.
Chose locally the second coordinate $w$ in such a way that the residue of
$\nabla$ is given by the formula

$$Res_{z=0}\nabla=
\left(\begin{array}{cc}
\beta-1&0\\
0&0\\
\end {array}\right)$$
In this case the connection $\nabla$ can be written as

\begin{equation}\label{lineconnection}
\nabla=d+\frac{dz}{z}Res_{z=0}\nabla+A,
\end{equation}
where $A$ is holomorphic. This proves the lemma.

\hfill $\square$

\begin{corollary}\label{torsionfree}
Suppose that $\Gamma(\Omega_S^1)=0$, i.e., there is no nontrivial
holomorphic 1-forms on $S$; then the torsion of $\nabla$ is
identically zero.
\end{corollary}
{\bf Proof.} The torsion of a meromorphic connection on the
tangent bundle to a complex manifold $X$ is a section of
$\Omega_X^2\otimes TX$.  In the case when $X$ is a two-dimensional
complex surface $S$ the bundle $\Omega_S^2\otimes TS$ is
two-dimensional and isomorphic to the bundle of holomorphic
$1$-forms $\Omega_S^1$. Thus on $S$ any holomorphic section of
$\Omega_S^2\otimes TS$ is identically zero.

\hfill $\square$

\subsection{Formulas for connections on $(1,1)$ cones}

In the following proposition we describe the connection
of a $PK$ metric in a neighborhood of a singular point of type
$(1,1)$.

\begin{proposition} \label{conformulas}
Consider a $PK$-cone of type $(1,1,\alpha)$
with linear coordinates $z,w$. For $i=1,...,n$ let $l_1=0,...,l_n=0$
be the equations of the singular lines of the cone.
Let $2\pi\beta_i$ be the conical angle at $l_i=0$. Then
the following holds:

1) The residue $\mathrm{Res}_{l_i}\nabla$ of $\nabla$ at $l_i$ is given
by a constant matrix-valued function $A_i$ and the connection
$\nabla$ is given by the following formula:

\begin{equation}\label{constantres}
\nabla=d+A=d+\sum_{i=1}^n A_i\frac{dl_i}{l_i}
\end{equation}

2) The residues $A_i$ satisfy the following equations:

\begin{equation}\label{3equations}
a)\; \sum_{i=1}^n A_i=(\alpha-1)\mathrm{Id};\quad b)\; 
\mathrm{tr}(A_i)=\beta_i-1;\quad
c)\; \{l_i=0\}=\mathrm{ker}(A_i)
\end{equation}

\end{proposition}

{\bf Proof of \ref{conformulas}.1) }

We prove first that the residue of $\nabla$ is constant
at any $l_i$. Indeed, the action of $\mathbb C^*$ on the cone
changes the $PK$ metric by a scalar factor, thus this action
preserves the connection $\nabla$. For any $c\in \mathbb C^*$ we have:
$$A(cz,cw)=A(z,w),$$
i.e.,  the residue ${\rm Res}_{l_i}(\nabla)$ is constant on $l_i$.

Consider now the following connection $\nabla'$ on $\mathbb
{C}^2$:
$$\nabla'=d+\sum_i{\rm Res}_{l_i}(\nabla)\cdot \frac{dl_i}{l_i}$$

We claim that $\nabla'=\nabla$. Indeed,  the matrix-valued $1$-form
$\nabla'-\nabla$ has no poles at the lines $l_i$ thus it is
holomorphic on $\mathbb{C}^2$. Moreover this 1-form is preserved
by the $\mathbb {C}^*$ action, i.e., it is identically zero.

\hfill $\square$

Next lemma is essential for the proof of \ref {conformulas}.2).
\begin{lemma}For the Euler field $e$ of the $PK$-cone we have:
$$\nabla_e e=e$$
\end{lemma}

{\bf Proof.} It is sufficient to check this identity for the flat
$\mathbb{C}^2$.

\hfill $\square$

{\bf Proof of \ref{conformulas}.2).}

Let us prove (\ref{3equations}a).
The Euler field $e$ on the $PK$-cone is given by

$$e=\frac{1}{\alpha}\left(z\frac{\partial}{\partial z}
+w\frac{\partial}{\partial w}\right)$$ We have

$$\nabla_e e=de(e)+\sum_{i=1}^n A_i(e)\frac{dl_i}{l_i}(e)=
\frac{1}{\alpha}e+\sum_{i=1}^n\frac{1}{\alpha}A_i(e)=e$$

This means that $\sum_{i=1}^n A_i(e)=(\alpha-1)e$, i.e.,
$A=(\alpha-1)\mathrm{Id}$.

Statements (\ref{3equations}b,c) are proven in Lemma \ref{half1}.

\hfill $\square$

\begin{remark} For $i=1,...,n$ let $l_i$ be lines in $\mathbb C^2$
containing the origin and let $\beta_i$ be complex numbers.
Then the space of matrixes $A_i$ satisfying 
(\ref{3equations}) has dimension $(n-3)$. 
Next formula gives the unique connection 
for $n=3$ with poles at the lines $z=0$, $w=0$, $z+w=0$

$$\nabla=d+\left(\begin{array}{cc}
(\beta_1-1)\frac{dz}{z}+\frac{\beta_2+\beta_3-\beta_1-1}{2}\frac{dz+dw}{z+w}
&\frac{\beta_2+\beta_3-\beta_1-1}{2}(\frac{dz+dw}{z+w}-\frac{dw}{w})\\
\frac{\beta_1+\beta_3-\beta_2-1}{2}(-\frac{dz}{z}+\frac{dz+dw}{z+w})&
\frac{\beta_1+\beta_3-\beta_2-1}{2}\frac{dz+dw}{z+w}+(\beta_2-1)\frac{dw}{w}\\
\end {array}\right)$$
\end{remark}

\begin{proposition}\label{complexaffine} 
Any connection $\nabla$ on $\mathbb C^2$ given by
formula (\ref{constantres}) with the matrices
$A_i$ satisfying (\ref{3equations}) is flat,
torsion-free and thus it defines a singular affine structure
on $\mathbb C^2$.
\end{proposition}

{\bf Proof.} Since  $dA=0$ the curvature of $\nabla$ is given by
$$dA+A\wedge A=A\wedge A$$
We need to prove that $A\wedge A=0$. Let us write
$$A_i=\left(\begin{array}{cc}
a_i&b_i\\
c_i&d_i\\
\end {array}\right),
$$
then the equation $A\wedge A=0$ is equivalent to the following system

$$\sum_i c_i\frac{dl_i}{l_i}\wedge \sum_i b_i\frac{dl_i}{l_i}=0 $$

$$\sum_i b_i\frac{dl_i}{l_i}\wedge \sum_i (a_i- d_i)\frac{dl_i}{l_i}=0$$

$$\sum_i c_i\frac{dl_i}{l_i}\wedge \sum_i (a_i- d_i)\frac{dl_i}{l_i}=0$$

For the first equation we have
$$\sum_i c_i\frac{dl_i}{l_i}\wedge \sum_i b_i\frac{dl_i}{l_i}=
d\mathrm{log} (\prod_i l_i^{c_i})\wedge d\mathrm{log}(\prod_i l_i^{b_i})$$

Function $f_1=\mathrm{log}(\prod_i l_i^{c_i})$
and $f_2=\mathrm{log}(\prod_i l_i^{b_i})$
are homogeneous of degree  $0$  on $\mathbb{C}^2$
(since $\sum c_i=\sum b_i=0$  by (\ref{3equations}a)).
It follows that $df_1\wedge df_2=0$. The next two equations are
completely analogous.

Now we show that $A$ is torsion free. We have

$$T\left(\frac{\partial}{\partial z},\frac{\partial}{\partial w}\right)=
\sum_{i=1}^n\frac{dl_i}{l_i}\left(\frac{\partial}{\partial
z}\right) A_i\left(\frac{\partial}{\partial w}\right)-
\frac{dl_i}{l_i}\left(\frac{\partial}{\partial w}\right)
A_i\left(\frac{\partial}{\partial z}\right)=$$
$$\sum_{i=1}^n\frac{1}{l_i}A_i\left(\frac{\partial l_i}{\partial z}
\frac{\partial }{\partial w}-\frac{\partial l_i}{\partial w}
\frac{\partial }{\partial z}\right)=0$$

The last equality holds by (\ref{3equations}c).

\hfill $\square$

{\bf Connection on the blow-up.}

\begin{lemma} \label{conblow}
Let $\nabla$ be a connection on $\mathbb C^2$ given by
formula (\ref{constantres}) with the matrices
$A_i$ satisfying (\ref{3equations}). Consider the blow-up of
$\mathbb C^2$ at $(0,0)$, take the pull-back of the tangent
bundle of $\mathbb C^2$ to the blow-up, and consider on it
the pull-back of $\nabla$. The obtained connection is
logarithmic and its residue at the exceptional curve is
$(\alpha-1)Id$.
\end{lemma}

{\bf Proof.} Let us introduce coordinates $u,v$ on the blow-up
$u=\frac{z}{w}, v=w$. The exceptional line is given by $v=0$.
Let $l_i=s_iz+t_iw=s_iuv+t_iv$. Then the pull-back connection
is given by
$$\nabla=d+\sum_{i=1}^n A_i d{\rm  (log} l_i)=
d+\sum_{i=1}^n A_i (d{\rm (log} v)+d{\rm (log}(s_iu+t_i))=$$
$$d+(\alpha-1){\rm Id}\cdot d({\rm log}v)
+\sum_{i=1}^n A_i\cdot d{\rm log}(s_iu+t_i))$$
Here we use Equation (\ref{3equations}a). 
This proves the lemma.

\hfill $\square$

\subsection{Unitary flat logarithmic torsion free connection 
$\mapsto$ $PK$ metric}
In this subsection we consider only arrangements $(S,\Gamma_j)$ that
satisfy the property that $\Gamma_j$ are smooth and transversal.
We give a sufficient criterion for such an arrangement to be the 
singular locus of a $PK$ metric in terms of an {\it adapted} connection.
Singularities of such an arrangement are normal crossings or 
singularities of type $(1,1)$.

\begin{definition} \label{parabconditions}
Let $\nabla$ be a connection partially 
adapted to $(S,\Gamma_j,\beta_j)$. Suppose that $\Gamma_j$ are 
smooth and intersect transversally. 
Consider the blow up $\pi:\widetilde S\to S$
at all points $x_i$ of multiplicity at least $3$. 
$\nabla$ is {\it adapted} to  $(\Gamma_j,\beta_j)$ if
the pull-back connection $\pi^*\nabla$ on $\pi^*TS$ is logarithmic
on $\widetilde S$ and its residue at the exceptional 
curve over $x_i$ equals $\sum_{j}d_{ij}(\beta_i-1)Id$.

\end{definition}

\begin{theorem} \label{connectionPK}
Let $(S,\Gamma_j, \beta_j, x_i)$ be a weighted 
arrangements of curves and $\nabla$ be an adapted 
flat, unitary, torsion free connection on $TS$. Suppose that 
$0<\beta_j<1$ and for every $x_i$, $\sum_i d_{ij}(\beta_i-1)>-2$.
Then the unitary metric on $TS$ corresponding to  $\nabla$ is a $PK$ metric.
\end{theorem}

{\bf Proof.} Note first that since $\nabla$ is unitary, flat, and 
torsion-free, the metric $g$ corresponding to $\nabla$
is flat on the complement to the curves $\Gamma_j$. So to prove that 
$g$ extends to a $PK$ metric on $S$ it will be sufficient to show 
the following $3$ properties of $g$.

a) For any  smooth point  of $\Gamma$  there is a 
neighborhood $U$ with $U\setminus \Gamma$
isometric to the direct product of a flat punctured $2$-cone 
and a flat disc. 

b) For any  double point  of $\Gamma$  there is a 
neighborhood $U$ with $U\setminus \Gamma$
isometric to the direct product of two flat punctured $2$-cones.

c) For any point of $\Gamma$ of multiplicity greater than $2$
a neighborhood of the point is isometric to a 
$2_{\mathbb C}$-dimensional $PK$ cone.

The proofs of a) and b) are similar so we will prove only b) and c).

{\bf Proof of b).} 
Introduce coordinates $z,w$ in a neighborhood of the double point 
such that $\Gamma$ is locally given by $zw=0$. According to 
Proposition \ref{malgrange} there exist two sections $s_1$ and $s_2$
such that the connection $\nabla$ is given by 
$$\nabla =d+B_1\frac{dz}{z}+B_2\frac{dw}{w}$$
where $B_1$ and $B_2$ are constant and commuting. Making an additional
linear change in $s_1,s_2$  we can 
suppose that $B_1(s_1)=0$ and $B_2(s_2)=0$. It is clear then that 
the sub-bundles of the tangent bundle generated by $s_1$ and $s_2$
are invariant under $\nabla$ and moreover the vector  field $s_1$
is tangent to the line $w=0$ and $s_2$ is tangent to the line $z=0$.   
Since $\nabla$ is torsion free the integral curves of $s_1$ and 
$s_2$ are flat cones. So we can deduce that 
locally the neighborhood of $(0,0)$ is 
the direct product of two $2$-cones.

\hfill $\square$

{\bf Proof of c).} Let $0$ be a point of $\Gamma$ of 
multiplicity greater than $2$. From a) it follows that
the metric  $g$  extends continuously to any punctured curve
$\Gamma_j\setminus 0$, thus we obtain a polyhedral K\"ahler metric on 
$U\setminus0$. It is necessary to show further that $g$ extends to 
$0$ and the resulting metric on $U$ is a polyhedral K\"ahler metric.
To show this it is sufficient to construct
the action of $\mathbb C^*$ on $U\setminus 0$ by dilatations, i.e., 
$\mathbb R^*$  must act  by dilatations of the metric 
and $S^1$ must act by isometries.

{\bf Construction of the $\mathbb C^*$ action on $U\setminus 0$.}

Consider the holonomy
representation of the group $\pi_1(U\setminus \Gamma, x)$ in the
group $U(2,\mathbb C)\ltimes \mathbb C^2$ of unitary affine
transformations of $\mathbb C^2$.
$$\mathrm{Hol}_{\nabla}:\pi_1(U\setminus \Gamma, x)
\to U(2,\mathbb C)\ltimes \mathbb C^2$$
The linear part of this representation is
denoted by $\mathrm {hol}_{\nabla}$. Denote by $c$ the generator of
the center of $\pi_1(U\setminus \Gamma, x)$ corresponding to an
anti-clockwise path around $0$ in a complex line through $0$.  
From condition 2) of 
Definition \ref{parabconditions} it follows that
$$\mathrm {hol}_\nabla(c)=\mathrm {exp}(2\pi i\sum_j(\beta_j-1))\mathrm {Id}$$

We see that the affine transformation ${\mathrm Hol}_\nabla(c)$
is a complex rotation of $\mathbb C^2$ around some point  $y$
on angle $2\pi\sum_j(\beta_j-1)$.

Since any element $h$ of $\pi_1(U\setminus \Gamma, x)$ commutes
with $c$, the affine transformation $\mathrm {Hol}_\nabla(h)$
fixes the point $y$. It follows that the representation $\mathrm
{Hol}_\nabla(c)\pi_1(U\setminus \Gamma, x)$ commutes with the
affine action of $\mathbb C^*$ on $\mathbb C^2$, given by complex
dilatations that fixe $y$. We deduce that  
the action of $\mathbb C^*$ can  be pulled back 
to the action of $\mathbb C^*$ on $U\setminus 0$. 

\hfill $\square$

\begin{remark} The condition $\beta_j<1$ in this theorem can 
be replaced by the  condition $\beta_j\notin \mathbb Z_+$
but we don't prove this here. 
\end{remark}

\section{Topological relations on the singular locus.}

In this section we prove Theorem \ref{Thm:chernclasses}.
We  use  a formula of  Ohtsuki [Oh].

{\bf Theorem [Oh].} Let $S$ be a surface and $E$ be a holomorphic
vector bundle  on it. Let $D=\bigcup_j D_j$  be a normal crossings
divisor on $S$ and $\nabla$ a  logarithmic connection on $E$
with poles at $D$. Denote by $y_k$ the double points  of $D$ and by
$D_{k_1}$ and $D_{k_2}$ the irreducible components of $D$ 
containing $y_k$. 
We also use the notation ${\rm R}_j={\rm Res}_{D_j}(\nabla)$.
The following identities hold:

$$c_1(E)=-\sum_j {\rm Tr}({\rm R}_j)D_j$$
$$c_2(E)=\sum_k ({\rm Det}({\rm R}_{k1}+{\rm R}_{k2})
-{\rm Det}({\rm R}_{k1})-{\rm Det}({\rm R}_{k2}))(y_k)
+
\sum_j {\rm Det}({\rm R}_j) D_j\cdot D_j$$

Note that the value of the first summand is defined only at $y_k$
but the function ${\rm Det(R}_j)$ is defined on the whole divisor $D_j$
and is constant on it.

\subsection{Proof of \ref{Thm:chernclasses}(2) and
\ref{Thm:chernclasses}(3) in the case of $(1,1)$ singularities}

Consider a $PK$ surface $S$ with the weighted arrangement
$(\Gamma_j,\beta_j; x_i,\alpha_i)$ such that all
multiple points of the arrangement are either normal
crossings or singularities of type $(1,1)$.
Let $S_b$ be the blow-up of $S$ 
at the points $x_i$  and let $\pi: S_b\to S$ be the blow-down.
Denote by $P_i$ the exceptional curve over $x_i$ and denote by $\widetilde
\Gamma_j$ the proper transform of $\Gamma_j$. Consider the
pull-back $\widetilde\nabla$ of the $PK$ connection $\nabla$ to
$\pi^*(TS)$. By Lemma \ref{conblow}
 $\widetilde\nabla$ is logarithmic with 
poles at the divisor $\bigcup_i P_i\bigcup_j\widetilde\Gamma_j$
(further we call this divisor by $D$), 
and the residue of $\widetilde\nabla$
at $P_i$ is equal to $(\alpha_i-1)Id$.

\begin{lemma}\label{pulback}
$$\pi^*\sum_j (\beta_j-1) \Gamma_j=
\sum_i 2(\alpha_i-1)P_i+\sum_j(\beta_j-1)\widetilde\Gamma_j$$
\end{lemma}
{\bf Proof.}
$$\pi^*(\sum_j (\beta_j-1)\Gamma_j)=
\sum_j(\beta_j-1)(\widetilde\Gamma_j +\sum_i d_{ij}P_i)=$$
$$=\sum_j(\beta_j-1)\widetilde\Gamma_j+\sum_i 2(\alpha_i-1)P_i$$
the second equality follows from  \ref{Thm:alpha=}.

\hfill $\square$

{\bf Proof of \ref{Thm:chernclasses}(2)}

Using first Theorem [Oh] and then Lemma \ref{pulback} we obtain:

$$c_1(\pi^*K_S)=-c_1(\pi^*TS)=\sum_i {\rm Tr(Res}_{P_i}\widetilde\nabla)
 \cdot P_i
+\sum_j {\rm Tr(Res}_{\widetilde\Gamma_j}\widetilde\nabla)
\cdot \widetilde\Gamma_j=$$

$$=\sum_i 2(\alpha_i-1)P_i+\sum_j(\beta_j-1)\widetilde\Gamma_j=
\pi^*\sum_j(\beta_j-1) \Gamma_j$$

It follows that $K_S=\sum_j(\beta_j-1) \Gamma_j$.

\hfill $\square$

{\bf Proof of \ref{Thm:chernclasses}(3)}

According to Theorem [Oh] the number $c_2(\pi^*TS)$ can be expressed as 
the sum of the contributions of the irreducible components of $D$ 
and the sum over their pairwise intersections. The first sum is the following:
$$\sum_i \mathrm{Det}(\mathrm{Res}_{P_i}\widetilde\nabla)P_i\cdot P_i+
\sum_j \mathrm{Det}(\mathrm{Res}_{\widetilde \Gamma_j}\widetilde\nabla)
\widetilde\Gamma_j\cdot\widetilde\Gamma_j=-\sum_i(\alpha_i-1)^2$$
Here we use
$\mathrm{Det}(\mathrm{Res}_{\widetilde \Gamma_j}\widetilde\nabla)=0$.

The sum of the contributions of the double points on $P_i$ is the following:
$$(\alpha_i-1)\sum_j d_{ij}(\beta_j-1)
=2(\alpha_i-1)^2$$

Any intersection of $\widetilde \Gamma_j$ with $\widetilde
\Gamma_k$ contributes $(\beta_j-1)(\beta_k-1)$ thus the
sum over all intersections of curves  $\widetilde \Gamma_j$ is
given by
$$\sum_{j>k}B_{jk}(\beta_j-1)(\beta_k-1)$$
Finally, taking the sum of all contribution we obtain

$$c_2(TS)=c_2(\pi^*(TS))=\sum_i(\alpha_i-1)^2+\sum_{j>k}
B_{jk}(\beta_j-1)(\beta_k-1)$$
\hfill $\square$

\subsection{Proof of \ref{Thm:chernclasses}(1)}

\begin{lemma} \label{lemma:[Gj:Gk]}
For any $j\ne k$ we have

\begin{equation}\label{GjGk}
\Gamma_j\cdot \Gamma_k=
B_{jk}+\sum_i (p_iq_i)(d_{ij}d_{ik})
\end{equation}
\end{lemma}

{\bf Proof.} This formula expresses the intersection index of  
$\Gamma_j$  and $\Gamma_k$ as a sum of local multiplicities of their 
intersections. By definition
$B_{jk}$ is the number of transversal 
intersections of $\Gamma_j$ and $\Gamma_k$. The local 
multiplicity of the  intersection of $\Gamma_j$ and $\Gamma_k$ at $x_i$ 
equals $(p_iq_i)(d_{ij}d_{ik})$. Indeed, 
the local multiplicity of the intersection of the curves $cz^q=w^p$  
and $z^q=w^p$ at $0$ equals $pq$ if $1\ne c\ne 0$; the 
local multiplicities of intersection of $z^q=w^p$ with lines $z=0$ 
and $w=0$ are equal to $p$ and $q$ correspondingly. 
Now everything follows from the definition of $d_j$ and $d_k$.

\hfill $\square$

{\bf Proof of \ref{Thm:chernclasses}(1)}

The Gauss-Bonnet theorem for flat surfaces with 
conical singularities implies

$$2g(\Gamma_j)-2=\sum_{k\ne j}B_{jk}(\beta_k-1)+
\sum_i(d_{ij}\alpha_i-\tilde d_{ij})$$

This formula expresses the Euler characteristics of $\Gamma_j$ as
the sum of the defects of the conical points of $\Gamma_j$. The first sum
contains the contribution of the normal crossings of $\Gamma_j$ and the
second sum contains the contribution of the singularities $x_i$.

Now, using Theorem \ref{Thm:alpha=}, we obtain  the following
expression for the right term of the previous equation:

$$\sum_{k\ne j}B_{jk}(\beta_k-1)+
\sum_i\left(d_{ij}\frac{p_i q_i}{2}\sum_k d_{ik}(\beta_k-1)
+d_{ij}\frac{p_i+q_i}{2}-\tilde d_{ij}\right)=$$

using Lemma \ref{lemma:[Gj:Gk]}
$$=\sum_{k\ne j}B_{jk}\frac{\beta_k-1}{2}+
\sum_k\Gamma_j\cdot \Gamma_k\frac{\beta_k-1}{2}-
\Gamma_j\cdot\Gamma_j\frac{\beta_j-1}{2}+$$

$$+\sum_i\frac{p_i q_i}{2}(d_{ij})^2(\beta_j-1)+
\sum_i\left(d_{ij}\frac{p_i+q_i}{2}-\tilde d_{ij}\right)=$$

using relation \ref{Thm:chernclasses}(2)
$$=\sum_k B_{jk}\frac{\beta_k-1}{2}+\frac{K_S\cdot\Gamma_j}{2}+
\sum_i\left(d_{ij}\frac{p_i+q_i}{2}-\tilde d_{ij}\right)$$

This proves \ref{Thm:chernclasses}(1).

\hfill $\square$

\subsection{Line arrangements in $\mathbb{C}P^2$}

In this subsection we consider weighted arrangements
of lines on $\mathbb{C}P^2$ that satisfy  
Equations \ref{Thm:chernclasses}(1)-(3).
Let  $(L_1,\beta_1;...;L_n,\beta_n)$
be such an arrangement. Any singularity of this arrangement is
either a normal crossing or a $(1,1)$-type singularity. So 
Equations \ref{Thm:chernclasses}(1) and \ref{Thm:chernclasses}(2)
simplify and take the form:

$$\sum_k B_{jk}(\beta_k-1)=-1,\;\;\;\;\sum_k(\beta_k-1)=-3$$

Here by definition the number  $B_{ij}$ $(i\ne j)$ is equal to $1$
if the point of the intersection of $L_i$ and $L_j$ is a double point
(i.e., other lines of the arrangement don't contain it) and
$B_{ij}=0$ otherwise. The number $(B_{jj}+1)$  is equal to the
number of points of multiplicity at least $3$ on the line $L_j$.

{\bf Symmetric case.} Consider the most symmetric case 
when all angles $\beta_j$
are equal. Then we have

$$\beta_j-1=-\frac{3}{n},\;\; \sum_k B_{jk}=\frac{n}{3}$$

The number $\sum_k B_{jk}+1$ is equal to the number of all
intersections of $L_j$ with other lines. Thus we obtain the
following condition:

{\it The arrangement contains $3m$ lines and any line intersects
the other lines exactly at $m+1$ points.}

One can show that such arrangements satisfy 
as well Equation \ref{Thm:chernclasses}(3).
These arrangements were considered first by Hirzebruch in [Hir] and 
we recall several examples. 

{\it 1) 3 lines.}  A generic configuration of $3$ lines on $\mathbb{C}P^2$.

{\it 2) 6 lines.}  The configuration of $6$ lines
$x-y=0$,  $x\pm z=0$, $y\pm z=0$ $z=0$.

{\it 3) $3(m+1)$ lines, $m>1$.} Consider the ramified covering of
$\mathbb{C}P^2$ by itself given by $(x:y:z)\to(x^m:y^m:z^m)$. The
preimage of the configuration of $6$ lines is an arrangement of
$3(m+1)$ lines and any line has $m+2$ intersections with the other
lines.

{\it 4) Hesse arrangement.} Consider a nonsingular cubic in
$\mathbb{C}P^2$. It has $9$ points of inflections. There exist $12$
lines in $\mathbb{C}P^2$ that intersect the cubic exactly at
the points of inflections.

{\bf Criterium.} There exists one 
geometric condition that often permits 
to decide quickly that a given line arrangement is not $PK$.

\begin{proposition} \label{doublecrit} Every $PK$ 
arrangement $(L_j,\beta_j)$ in $\mathbb CP^2$ satisfy the following 
criterium. For every $L_j$ the exist a point in 
$\mathbb CP^2$ that belongs to all lines $L_k$ such that $L_k\cap L_j$ 
is a double point of the arrangement.
\end{proposition}

{\bf Proof.} Consider the sub-bundle of 
$T\mathbb CP^2|_{L_j}$ of directions orthogonal to $L_j$ with 
respect to $PK$ metric (these directions are 
eigenvectors of the residue map $\rm{Res}_{L_j}\nabla$).
This sub-bundle  is holomorphic outside of 
the multiple points of the arrangement and it 
extends holomorphically on the whole line $L_j$. Indeed, at
double points the $PK$ metric is a direct product of two $1$-cones, and 
at every point of multiplicity more than $2$ the eigenvectors of 
 $\rm{Res}_{L_j}\nabla$ are constant in the local linear coordinates
(Proposition \ref{conformulas}). 
The defined sub-bundle of $T\mathbb CP^2|_{L_j}$  is transversal
to $L_j$ and so it has degree $1$. It follows 
that there exists a point $y_j$  
in $\mathbb CP^2$  such that this sub-bundle of  
$T\mathbb CP^2|_{L_j}$  is given 
by directions tangent to the lines through $y_j$. This proves the proposition.

\hfill $\square$

\subsection{Limit $PK$ arrangements with a cusp}

\begin{definition} An arrangement of lines on $\mathbb RP^2$
is called {\it simplicial} if it subdivides $\mathbb RP^2$ in 
triangles.  
\end{definition}

Simplicial arrangements often occur as solutions to some 
extreme (combinatorial) problems [G2] and so it is not very 
surprising that some of simplicial arrangements are $PK$
($PK$ arrangement give extremum 
for the Bogomolov-Gieseker Inequality (\ref{cp2gieseker}) 
Theorem \ref{th:existence}).
For the moment $3$ infinite family of simplicial
arrangements and $91$ sporadic examples are know; $90$ sporadic   
examples are listed in [G1] and one additional in [G2].
It will be interesting to find out what sporadic arrangements 
are $PK$, the criterium from Proposition \ref{doublecrit}  
rules out some of them.

For three infinite series of arrangements from [G1] 
there exists a system of weights such that  all {\it equalities}
of Theorems \ref{Thm:chernclasses} and \ref{th:existence} are satisfied.
The first series is a union of a pencil of $n$ lines and 
a line that does not belong to the pencil. For this series 
the weight of the line should be $0$ and the weights of 
the lines from the pencil can be arbitrary (but we impose of course $\sum_j(\beta_j-1)=-3$).
Second series is called $R_{2k}$ and consists of the lines formed by extending
the sides of a regular $k$-gon together with an additional $k$ lines
formed be the axes of symmetry of the $k$-gon. We associate to 
the axes weight $\frac{k-1}{k}$ and to the sides weight $\frac{2k-1}{2k}$
(the choice of wights is unique for $k>3$). The third series $R_{4k+1}$
is the union of $R_{4k}$ with the line at infinity of weight $1$.
We see that for the point of the highest multiplicity of these 
arrangements the inequality $\sum_{ij}d_{ij}(\beta_j-1)>-2$ 
does not hold strictly but instead of this the equality holds. 
So these arrangements are not $PK$.
We conjecture instead that there is a different geometric structure 
related to these arrangements.

\begin{definition} A weighted arrangement $(L_j,\beta_j)$ 
in $\mathbb CP^2$ is called a {\it limit $PK$ 
arrangement  with a cusp} if it satisfies all conditions 
of Theorem \ref{th:existence} apart from one inequality. 
We impose that there is 
a multiple point $x$ of the arrangement called  the {\it cusp}
such that $\sum_{j|x\in L_j} (\beta_j-1)=-2$. 
\end{definition}

Arrangements $R_{2k}$ satisfy this definition, $R_{4k+1}$
and a pencil of lines plus a line formally don't satisfy 
(because the first arrangement has a line of  weight $1$ and 
the second has a line of  weight $0$). But the following should hold 
for all $3$  series.

\begin{conjecture}\label{pencil} For every limit $PK$
arrangement with a cusp
there exists a flat torsion-free connection on $\mathbb CP^2$ with 
holonomy in the upper triangular subgroup of 
$SL(2,\mathbb C)$ and with the poles of residues $(0,\beta_j-1)$
at the lines $L_j$. This connection should preserve the sub-line bundle
of $T\mathbb CP^2$ tangent to the pencil of lines through the cusp
of the arrangement. 
\end{conjecture}

For a pencil of lines plus a line the connection on $\mathbb CP^2$ 
should be given by a formula from Proposition \ref{complexaffine}.
Here $\mathbb CP^2$ is the completion of $\mathbb C^2$, 
the line at infinity belongs to the arrangement 
and has weight $\beta=0$.

If Conjecture \ref{pencil} holds it should be possible 
to deduce that limit $PK$ arrangements with a cusp satisfy the 
following restrictive properties. 

{\bf Conjectural properties.} 
For every multiple point $y$ 
of the arrangement the line $[x,y]$ belongs to the arrangement
(here $x$ is the cusp). If $mult(y)\ge3$ then 
sum of the defects of the 
lines  $L_j$ that contain $y$ but don't contain $x$  
equals the defect of the line $[x,y]$.

It will be very interesting to classify all weighted arrangements that satisfy
the two conjectural properties (they hold for $R_{2k}$ and $R_{4k+1}$). 
This will help to classify 
{\it non-rigid} $PK$ arrangements for which the admissible 
collection of weights $\beta_j$ have moduli. Of course the weights $\beta_j$
belong to a certain open polyhedron and the weights corresponding to
limit $PK$ arrangements can appear on the boundary 
of the polyhedron.







\section{Parabolic bundles and Kobayashi-Hitchin correspondence}

The goal of this section is to recall the notion of parabolic 
bundles and to formulate in a handy way several results from [M2] 
that we use in the proof of Theorem \ref{th:existence}.
In particular we formulate the parabolic version of Kobayashi-Hitchin 
correspondence from [M2]. A  systematic and thorough treatment of
parabolic bundles  can be found in [IS]. Parabolic Chern character
is also defined in [IS] (but we will not use it here).
We adopt partially the notations of these articles.

We will discuss only parabolic bundles on complex surfaces.
A good reference for usual two-dimensional bundles on surfaces is [F].

\begin{definition}
Let $X$ be a complex surface
and $D$ be a simple normal crossing divisor
with the irreducible decomposition $D=\bigcup_{i\in S} D_i$.
A {\it parabolic bundle} $E_*$ on $X$ is given by
a bundle $E$ with a collection of increasing filtrations by
sub-sheaves
$F_a^i$, indexed by $i\in S$, $a\in ]0,1]$ and
satisfying the following properties:

1) Every  sub-sheaf  $F_a^i$ is  locally free.

2) $ E(-D_i)\subset F^i_a(S)$ for any $a\in ]0,1]$.

3) The sets $\{a| F_a^i(E)/F^i_{<a}(E)\ne 0\}$
are finite for any $i$ in $S$.
\end{definition}

\begin{remark}
A parabolic bundle $E_*$ on a complex surface $X$ with 
a simple normal crossing divisor $D$ induces natural filtrations on 
the restrictions $E|_{D_i}$ by their vector sub-bundles. 
These filtrations are indexed by 
$a\in ]0,1]$ and 
defined by the formula:
$$F_a^i/E(-D)\subset E|_{D_i}$$
Parabolic structure can be reconstructed from these filtrations (see [IS]).

\end{remark}

\begin{definition} \label{parch1}

The {\it parabolic first Chern class} of a
parabolic bundle $E_*$ is given by the following formula:

\begin{equation}
{\rm par\textrm{-}ch_1}(E_*)={\rm ch_1}(E)-
\sum_i\sum_{a_i}a_i\cdot {\rm rank}_{D_i}(F^i_{a_i}/ F^i_{<a_i})\cdot [D_i]
\end{equation}

Let $L$ be an ample line bundle on $S$. Then the {\it parabolic degree}
of $E_*$ with respect to $L$ is given by
\begin{equation}
pardeg_L(E_*)=\int_S {\rm par\textrm{-}ch_1}(E_*) \cdot c_1(L)
\end{equation}

\end{definition}
{\bf The parabolic second Chern characteristic number}.

Let $(S,D)$ be a complex surface with a normal crossing
divisor $D=\cup_{i=1}^n D_i$. Let $E_{*}$ be a parabolic
vector bundle. We will recall now  the formula for the
parabolic second character  of $E_{*}$ given by
Mochizuki in [M2]. His formula works in much larger generality but
we need only the case of surfaces.

The parabolic second Chern character  of $E_*$ is
given as a sum  of $c_2(E)$, the contributions of the divisors
$D_i$, and the points of their intersections $D_i\cap D_j$.

To define the contributions of the points in $D_i\cap D_j$, for every
$a_i,a_j\in ]0,1]$,
consider the skyscraper sheaf $Gr^F_{(a_i,a_j)}$

$$F^i_{a_i}\cap F^j_{a_j}/
((F^i_{a_i}\cap F^j_{<a_j})\cup (F^i_{<a_i}\cap F^j_{a_j}))$$

This sheaf is  supported at the points in $D_i\cap D_j$,
and it is non-trivial only for finite set of $(a_i,a_j)$. Consider the
following sum

$$\nu(i,j)=\sum_{p\in D_i\cap D_j; a_i, a_j}a_i\cdot a_j
\cdot {\rm rank_p} (Gr^F_{(a_i,a_j)})$$

\begin{definition}\label{parch2def}
The second parabolic Chern character of the parabolic bundle $E_*$
is given by the formula
$$par\textrm{-}ch_2(E_*)=ch_2(E)-
\sum_{i;a_i}a_i\cdot c_1(F^i_{a_i}/F^i_{<a_i})+$$
$$+\sum_{i;a_i}\frac{1}{2}a_i^2{\rm rank}_{D_i}(F^i_{a_i}/F^i_{<a_i})
\cdot [D_i\cdot D_i]+\sum_{i,j}\nu(i,j)$$
\end{definition}

{\bf Stable bundles and Bogomolov-Gieseker inequality.}

Here again we consider  a surface $S$ with a parabolic
bundle $E_*$.  For any sub-sheaf
$V$ of $E$ the filtration on $E$ induces a structure of  a parabolic
sheaf on $V$. Recall that  a
sub-sheaf $V$ of $E$ is called {\it saturated} is the quotient
$E/V$ is torsion-free.

\begin{definition}Let $L$ be an ample bundle on $S$.  The bundle
$E_*$ is called {\it $\mu_L$-stable} (or {\it slope stable}) 
if for every saturated
sub-sheaf $V$ of $E$ it holds

$$pardeg_L V_* <pardeg_L E_*$$
\end{definition}

\begin{remark}
When $E$ is a rank $2$ bundle in order to
check its stability it is sufficient to consider only
{\it saturated locally free rank one sub-sheaves}  of $E$.
Following [F] we call such
sub-sheaves of $E$ {\it sub-line bundles}.
\end{remark}

The following inequality (called {\it Bogomolov-Gieseker} inequality) 
is proven in [M2] and was also proven in different terms in [Li].

\begin{theorem}\label{bginequality}
Let $E_*$ be a $\mu_L$-stable parabolic
 bundle on the surface
$S$. Then following inequality holds:
\begin{equation}
{\rm par\textrm{-}ch}_{2}(E_*)-
\frac{1}{2}{\rm par\textrm{-}ch}_1^2(E_*)\le 0
\end{equation}
\end{theorem}

\begin{definition} Let $(X,D,E_*)$ be a complex surface with 
a simple normal crossing divisor $D$ with a parabolic bundle $E_*$,
and let $\nabla$ be a unitary flat logarithmic connection on $E$ with
poles at $D$. We say that $\nabla$ is {\it compatible} with $E_*$
if the following conditions hold. 

1) For every $i$ and 
$a\in ]0,1]$ the sub-bundle $F_a^i/E(-D_i)$ of $E|_{D_i}$ is preserved 
by the residue map $Res_{D_i}(\nabla)$. 

2) The eigenvalue of $Res_{D_i}(\nabla)$
on the bundle $F_a^i/F_{<a}^i$ equals $-a$ (recall that this bundle 
is non-trivial only for finite number of values of $a$). 

3) The connection
induced by $\nabla$ on $F_a^i(X\setminus D_i)$ extends to a 
logarithmic connection on $F_a^i(X)$.

\end{definition}

Finally we can formulate  the version of parabolic Kobayashi-Hitchin 
correspondence that we use later.

\begin{theorem}\label{takuro} Let $(X,D,E_*)$ be a complex surface 
with  a simple normal crossing divisor $D$ and  
a parabolic bundle $E_*$. Suppose that $E_*$ is $\mu_L$-stable, 
has zero parabolic degree, and has zero second parabolic Chern number. 
Then there exists  a flat unitary logarithmic connection on 
$E$ compatible with $E_*$.
\end{theorem} 

This statement can be deduced from [M2] (see also [Li])  and [M1].
We explain this very briefly skipping all details. 
In  [M2] Mochizuki works  with parabolic Higgs bundles and proves  
Kobayashi-Hitchin correspondence for stable Higgs bundles with 
vanishing first and second Chern characters. 
The case  that we are interested in is the particular 
case when the Higgs field is zero. So by  
[M2] there exists a unitary flat metric on $E(X\setminus D)$ 
adapted to the parabolic structure. 
It follows further from [M1] that the flat unitary connection 
corresponding  to the flat metric constructed on $E(X\setminus D)$ 
extends to a logarithmic connection on $E$ and moreover this 
connection is compatible with the parabolic structure $E_*$.  

\section{Theorem of existence}
In this section we prove  Theorem \ref{th:existence}.

\subsection{A description of the proof}

Let  $(L_j,\beta_j)$ be a weighted  arrangement of lines in
$\mathbb{C}P^2$. Recall that by $x_i$ we denote the multiple
points of the arrangement of multiplicity at least $3$; 
$d_{ij}=1$ if $x_i$ belongs to $L_j$ and $d_{ij}=0$ otherwise.
Suppose that $(L_j,\beta_j)$ satisfies the three conditions 
of Theorem \ref{th:existence}.

To prove Theorem \ref{th:existence} we make the blow up
$\pi :S \to \mathbb CP^2$ of
$\mathbb CP^2$ at the points $x_i$
and consider on $S$ the pull back $E$ of the tangent bundle
$E=\pi^*T\mathbb {C}P^2$.
Using the  weights $\beta_j$ we define a  parabolic structure on
$E$. We get a parabolic bundle $E_*$ with
zero parabolic first Chern class and prove that $E_*$ is stable
for a certain  polarisation.
Inequality (\ref{cp2gieseker}) is just the Bogomolov-Gieseker 
inequality (Theorem \ref{bginequality}).

In the case when the second parabolic Chern class of $E_*$
equals $0$  using Theorem \ref{takuro} we prove that there exists 
a logarithmic flat unitary  connection on $T\mathbb CP^2$
and combining this with results of Section 4 we 
conclude the proof of Theorem \ref{th:existence}.

\subsection{The parabolic bundle on the blown up $\mathbb CP^2$}

Let  $(L_j,\beta_j)$ be an arrangement of lines in
$\mathbb{C}P^2$, satisfying the inequalities (\ref{stabin}). In
this subsection we construct the parabolic structure on a $2$-bundle
on the blowup of $\mathbb CP^2$ and we calculate its parabolic Chern
character. Let us fix some notations.

{\bf Notations}.

Denote  by $S$ the blow up  
of $\mathbb CP^2$ at the  points $x_1,...,x_k$ 
of multiplicity at least $3$ and let 
$\pi: S\to \mathbb{C}P^2$ be the corresponding projection map;

Denote by $E$ the pullback bundle $\pi^*T\mathbb{C}P^2$;

For $j\in \{1,...,n\}$  denote by $D_j$ the proper transform of $L_j$;

For $j\in \{n+1,...,n+k\}$ denote by $D_j$ the exceptional line
$P_{j-n}$;

$$D_j=P_{j-n},\;\; \pi(P_{j-n})=x_{j-n}$$

\medskip

{\bf The parabolic structure on the bundle $(S,E)$.}

For any $j\in\{1,...,n\}$ the restriction of $E$
to $D_j$ contains a distinguished rank $1$ sub-bundle -- the pullback
$\pi^*TL_j$ of the tangent bundle of $L_j$. Denote by $E_j$
the subsheaf  of $E$, generated by the sections that are contained
in $\pi^*TL_j$, being restricted to $D_j$. This subsheaf fits into
the following exact sequence:

$$0\to E_j\to E\to E|_{D_j}/\pi^*TL_j\to 0$$

Now for $1\le j \le n$ we put
$$F^j_a=E^j, \;\;{\rm for}\;\; 0<a<1-\beta_j, \;\;\;
F^j_a=E \;\;{\rm for}\;\; 1-\beta_j\le a\le 1  $$
And for $0<i\le k$
$$F^{n+i}_a=E(-D_{n+i}), \;\;{\rm for}\;\; 0<a<1-\alpha_i, \;\;\;
F^{n+i}_a=E \;\;{\rm for}\;\; 1-\alpha_i\le a\le 1  $$

\begin{proposition}\label{parchcalc}
The first and second parabolic Chern characters of $(E_*)$
are given by the following formulas:
\begin{equation}
{\rm par\textrm{-}ch}_1(E_*)=ch_1(E)-\sum_{j=1}^n(1-\beta_j)D_j-
 2\sum_{i=1}^k (1-\alpha_i)D_{n+i}
\end{equation}

\begin{equation}\label{parch2cp2}
{\rm par\textrm{-}ch}_2(E_*)=\frac{3}{2}
-\sum_{j=1}^n (1-\beta_j)-\sum_{j=1}^n\frac{1}{2}(1-\beta_j)^2B_{jj}
+\sum_{i=1}^k (1-\alpha_i)^2
\end{equation}
\end{proposition}

{\bf Proof.} The calculation of $\rm par\textrm{-}ch_1(E_*)$ is a
straight-forward application of Definition \ref{parch1}. For
$j\in \{1,...,n\}$ the quotient sheaf $F^j_a/F^j_{<a}$ has rank
$1$ for $a=1-\beta_j$ and is trivial otherwise.
For $j\in \{n+1,...,n+k\}$ the sheaf $F^j_a/F^j_{<a}$ has rank $2$
for $a=1-\alpha_j$ and is trivial otherwise.

In order to calculate $\rm par\textrm{-}ch_2(E_*)$ 
we need the following lemma.

\begin{lemma} 1) For  $1\le j_1<j_2\le n$ we have $\nu(j_1,j_2)=0$.

2) For   $1<i_1<i_2\le k$ we have $\nu(n+i_1,n+i_2)=0$.

3) For $1\le j\le n$,  $1\le i\le k$ we have
$\nu(j,n+i)=d_{ij}(1-\beta_j)(1-\alpha_j)$.
\end{lemma}

{\bf Proof.} 

1) The sheaf
$Gr^F_{(a_{j_1},a_{j_2})}$ can be nontrivial only  when
$a_{j_1}=1-\beta_{j_1}$ and  $a_{j_2}=1-\beta_{j_2}$. 
But in this case  by construction the sheaf is equal to $E/(E_{j_1}\cup E_{j_2})$.
At the same time $E_{j_1}\cup E_{j_2}=E$.

2) $Gr^F_{(a_{n+i_1},a_{n+i_2})}$  is trivial   because
$D_{n+i_1}\cap D_{n+i_2}=\varnothing$.

3) The sheaf  $Gr^F_{(a_{j},a_{n+i})}$ is nontrivial only when
$a_{j}=1-\beta_j$, $a_{n+i}=1-\alpha_i$. In this case it is equal to
$E/(E(-D_{n+i})\cup E_j)$. It is supported at the points  $D_{j}\cap D_{n+i}$
and has rank $1$ at each point.

\hfill $\square$

{\bf Calculation of $\rm par\textrm{-}ch_2(E_*)$.} According to 
Definition (\ref{parch2def}) we have:

$$par\textrm{-}ch_2(E_*)=\frac{3}{2}-
\sum_{j=1}^n (1-\beta_j)$$
$$-\sum_{j=1}^n\frac{1}{2}(1-\beta_j)^2B_{jj}-
\sum_{i=1}^k (1-\alpha_i)^2
+\sum_{i=1,j=1}^{k,n}d_{ij}(1-\beta_j)(1-\alpha_i)$$

Here we use the following facts:

1) For $j\in (1,...,n)$  we have $c_1(E/E_j)=1$, $[D_j\cdot D_j]=-B_{jj}$.

2) For $i\in (1,...,k)$ we have $E|D_{n+i}$ is trivial, $[D_{n+i}\cdot D_{n+i}]=-1$.

Finally, by Theorem \ref{Thm:alpha=} (use $p=q=1$) the 
last term is equal to $2\sum_{i=1}^k (1-\alpha_i)^2$.
This concludes the proof.

\hfill $\square$

\subsection{Theorem on stability and additional lemmas}

In this subsection we prove that the parabolic bundle $E_*$
constructed above is stable with respect
to an appropriate  polarisation on $S$.
Take $N\in \mathbb Z_+$,
$$N>max\{\frac{k}{min_i \alpha_i},\;
\frac{2k}{min_{j,k}(\beta_j+\beta_k)},\;
\frac{3k}{1-max_j \beta_j}\}$$
and define the following line bundle:

$$L_N=O_S(-\sum_{i=1}^k P_i)\otimes p^*O(N)$$
Note, that $L_N$ is ample since $N>k$.

\begin{theorem} \label{stabth}
The parabolic bundle $(E_*,S)$ is $\mu_{L_N}$-stable. 
\end{theorem}

Let us first give a plan of the proof. The parabolic 
degree of $E_*$ is zero, so we need to show that  the parabolic degree
of any saturated sub-line bundle $V$ of $E$ is negative.
Every line bundle on $S$ is of the form 
$O_S(\sum_i d_iP_i)\otimes \pi^*O(d)$. For sub-line bundles of $E$ we have $d\le 1$ and the constant $N$ is chosen in such a way that the 
degree of $V\subset E$ can be positive only for $d\ge 0$. 
Since the parabolic weights are in $[0,1]$, $pardeg V_*\le deg V$,
so we only need to consider line sub-bundles of $E$ with $d=1,\;0$.
In the case $d=1$ the pushdown $\pi_*V$ is contained in 
a sub-line bundle of $T\mathbb CP^2$ generated by sections 
tangent to a pencil of lines. We prove that $pardeg (V_*)$ is negative 
comparing the degree of $V$ with the parabolic contribution, 
coming from the behavior of the pencil of lines with respect 
to the line arrangement on $\mathbb CP^2$.

Let us introduce some notations.
For a point $x$ in $\mathbb CP^2$ denote by $L_x$ the sub-line bundle 
of $T\mathbb CP^2$ generated by the sections tangent to the pencil of 
lines containing $x$. For a section $v$ of $T\mathbb CP^2$ with 
isolated zeros denote by $L_v$ the sub-line bundle of $T\mathbb CP^2$
generated by $v$. The following lemma is standard, we omit the proof.

\begin{lemma} 

Sub-line bundles of $T\mathbb CP^2$ have degree at most $1$.
Every saturated sub-line bundle of $
T\mathbb CP^2$ of degree $1$ equals $L_x$
for some $x$. Every saturated sub-line bundle of $T\mathbb CP^2$ 
of degree $0$ equals $L_v$ for some $v$.

\end{lemma}

\begin{lemma} \label{deg10}
Let $L=O_S(\sum_{i=1}^k d_iP_i)\otimes\pi^*(O(d))$
be a saturated sub-line bundle of $E$. Then $d\le 1$. Suppose
$d=1$ or $d=0$.

1) If $d=1$ then $\pi_*(V)\subset L_x$ for some  
$x\in \mathbb CP^2$.

2) If $d=0$ then $\pi_*(V)\subset L_v$
for some vector filed $v$ with isolated zeros.

3) $\pi_*(V)$  coincides with $L_x$ or $L_v$ outside of the set $\{x_1,...,x_k\}$.
\end{lemma}

{\bf Proof.} Consider $(\pi_*V)^{\vee\vee}$, this is 
saturated sub-line bundle of $T\mathbb CP^2$. Its degree equals to $d$,
so $d\le 1$. In the case $d=1$ by the previous lemma $(\pi_*V)^{\vee\vee}$
is $L_x$ of some $x$, if $d=0$ it is $L_v$ for some $v$.
The sheaf $\pi_*V$ is a subsheaf of $(\pi_*V)^{\vee\vee}$ and it 
differs from it only at points $x_i$ for which $d_i<0$.

\hfill $\square$

\begin{lemma} \label{subbundles}

Let $V=O_S(\sum_{i=1}^k d_iP_i)\otimes\pi^*(O(d))$
be a sub-line bundle of  $E$. Then 
for any $i$  we have $d_i\le 2-d$. In particular
we have an upper bound  on degree of $V$:

\begin{equation}\label{degin}
{\rm deg}_{L_N}(V)=c_1(V)\cdot c_1(L_N)\le (2-d)k+Nd
\end{equation}

\end{lemma}

{\bf Proof.} 
Let us prove that for any $i$ it holds $d_i\le 2-d$.
For any line $P$ in  $\mathbb CP^2$ we have
$T\mathbb CP^2_{|P}\simeq O(1)\oplus O(2)$. Take a line $P$ that contains
a point $x_i$ and doesn't contain any point $x_j$ for $j\ne i$.
Let $P'$ be the
proper transform of $P$. Then  again $E_{| P'}\simeq O(1)\oplus O(2)$.
Since $\mathrm{Hom_{O_S}}(V,E)\ne 0$, there is a line  $P$ through
$x_i$ for  which $\mathrm{Hom_{O_{P'}}}(V_{|P'},E_{|P'})\ne 0$.
At the same time, by definition of $V$
$$V_{|P'}=O(d\pi^*H+d_iP_i)_{|P'}=O(d+d_i)$$
it follows
$$(d+d_i)\le 2$$
Now we conclude the proof
$${\rm deg_{L_N}}(O_S(\sum_i d_iP_i+d\pi^*H))=
\sum_id_i+Nd\le (2-d)k+Nd$$
 \hfill $\square$

\begin{lemma}\label{sumpencil}
Let $(L_j,\beta_j)$ be a weighted line arrangement satisfying 
the conditions (\ref{stabin}) of Theorem \ref{th:existence}.

1) For any point $x$ in $\mathbb CP^2$
the following inequality holds:

\begin{equation}\label{stab}
\sum_{j|x\notin L_j}(1-\beta_j) > 1+\frac{2k}{N}
\end {equation}

2) A holomorphic vector field $v$ with isolated zeros can be tangent
to at most $3$ lines of the arrangement.

\end{lemma}

{\bf Proof.} 1)
It is clear that the sum attains its maximum for a point
that is a multiple point of the arrangement.
Since $\sum_j(1-\beta_j)=3$, for a double point of the 
arrangement the sum in (\ref{stab}) 
is at least $1+min_{j,k}(\beta_j+\beta_k)$.
For a point  $x_i$  of multiplicity more than $2$ we have:
$$\sum_{j|x\notin L_j}(1-\beta_j)=3-\sum_j d_{ij}(1-\beta_j)=
3-2(1-\alpha_i) \ge 1+2\mathrm{min}_i(\alpha_i)$$

2) This is standard, the field $v$ must have $3$ zeros and
it must be tangent to the lines that join these zeros.

\hfill $\square$

\subsection{Proof of stability}

{\bf Proof of Theorem \ref{stabth}.}

By \ref{deg10} and \ref{subbundles} any sub-line bundle
$V$  of $E$  is of the form $O_S(\sum_{i=1}^k d_iP_i)\otimes\pi^*(O(d))$
with  $d\le 1$, and $deg_{L_N}(V)<0$  if $d<0$.
Since in  our situation the parabolic weights are
contained in $]0,1]$, we have an inequality
$pardeg_{L_N}(V_*)\le deg_{L_N} (V)$.
So it is necessary only to consider the cases
when $d=1$ and $d=0$.

To calculate the parabolic first
Chern class $parch_1(V_*)$ we need to find
for every $j\in \{1,...,n+k\}$ and $a\in ]0,1]$  the
rank of the following quotient sheaf supported on $D_j$:
$$(V\cap F^j_{a_j})/(V\cap F^j_{<a_j})$$

Consider the case $d=1$. Then according to Lemma \ref{deg10}
there exists $x\in \mathbb CP^2$ such that $\pi_*(V)\subset L_x$.

In the case $j\in\{1,...,n\}$
there are two possibilities.

If $x\in L_j$
then $V\subset F^j_{a_j}$ for all $0<a_j\le 1$ so the corresponding
quotient sheaf is always trivial. If $x\notin L_j$ then
the quotient sheaf is non-trivial for $a_i=1-\beta_j$ and has
rank one (this follows from \ref{deg10}).

In the case $j\in\{n+1,...,n+k\}$ the quotient sheaf
is nontrivial for $a_j=1-\alpha_{j-n}$ and has rank $1$.
This gives us the formula:

$$parch_1(V_*)=ch_1(V)-\sum_{j|x\notin L_j}(1-\beta_j)D_j-
\sum_{i=1}^n(1-\alpha_i)D_{n+i} $$

We have the following sequence of inequalities
$$pardeg_{L_N}(V_*)=deg_{L_N}(V)-
c_1(L_N)\cdot (\sum_{j|x\notin L_j}(1-\beta_j)D_j+
\sum_{i=1}^n(1-\alpha_i)D_{n+i})\le$$
using Inequality (\ref{degin}) and the equality
$2(1-\alpha_i)= D_{n+i}\cdot \sum_j(1-\beta_j)D_j $
$$\le N+k-N\sum_{j|x\notin L_j}(1-\beta_j)+\sum_i 2(1-\alpha_i)-\sum_i(1-\alpha_i)<$$
using $0<\alpha_i<1$ and Inequality (\ref{stab})
$$< N+2k-N(1+\frac{2k}{N})=0$$

The case $d=0$ is analogous. By Lemma
\ref{deg10} there exists a vector-field $v$ with
isolated zeros such that  $\pi_*(V)\subset L_v$.
The arrangement contains more than $3$ lines
so there is at least one line that is not tangent to $v$.
Making the same calculation in the case $d=1$ and using 
$N>\frac{3k}{1-max_j \beta_j}$ we get

$$pardeg_{L_N}(V_*)\le 2k-N(1-max_j(\beta_j))+\sum_i(1-\alpha_j)<0$$

\hfill $\square$

{\bf Example.}
Consider the arrangement of $6$ lines that pass through
$4$ generic points $x_1,...,x_4$. For every $\beta\in]0,1[$
we can associate the weight $\beta$ to the lines $x_1x_i$
and $1-\beta$ to the lines $x_ix_j$, $i,j>1$. We get a stable
parabolic bundle on $\mathbb CP^2$ blown up at $x_1,...,x_4$. 
When $\beta$ tends to $0$ the parabolic degree of the sheaf corresponding
to $L_{x_1}$ tends to zero and as a result for $\beta=0$ we don't
get a $PK$ metric on $\mathbb CP^2$.

\subsection{Proof of Theorem \ref{th:existence} and an application 
of the Bogomolov Gieseker inequality}

{\bf Proof of Theorem \ref{th:existence}.} 
Let us sum up what we have done. We started with a weighted arrangement
$(L_j,\beta_j)$ that satisfies the conditions of Theorem \ref{th:existence}.
We introduced a structure of parabolic bundle $E_*$ on the 
pullback $E$ of the tangent bundle $T\mathbb CP^2$ to
the blow up of $\mathbb CP^2$ at the multiple points of the arrangement.
We proved that $E_*$  is stable (Theorem \ref{stabth})
and has zero first parabolic Chern class (Proposition \ref{parchcalc}).
So inequality (\ref{cp2gieseker}) follows from the calculation of 
second parabolic Chern number (Proposition \ref{parchcalc}) and Bogomolov-Gieseker
inequality (Theorem \ref{bginequality}). This proves the first part 
of the theorem.

If the equality is attained in (\ref{cp2gieseker}) then the   
second parabolic Chern number of $E_*$  vanishes, so we can use
parabolic Kobayashi-Hitchin correspondence 
(Theorem \ref{takuro}).
Namely, there exists a unitary flat logarithmic connection on $E$
compatible with the parabolic structure of $E_*$.
By Corollary \ref{torsionfree} the constructed connection is torsion-free.
Finally, using Theorem \ref{connectionPK} we conclude 
that the corresponding flat connection on 
$T\mathbb CP^2$ defines a $PK$ metric.

\hfill$\square$

Let us give one corollary of Theorem \ref{th:existence}.
For a multiple point $x_i$ of a line arrangement $L_1,...,L_N$
denote by $\mu_i$ the number of lines through $x_i$.

\begin{corollary} Suppose that the multiplicity of every point 
of the arrangement $(L_1,...,L_N)$ is less than $\frac{2N}{3}$.
Then the following inequality holds:
$$\sum_i\mu_i\ge \frac{N^2}{3}+N $$
In the case of equality $N$ is divisible by $3$, every 
line intersects over lines in $\frac{N}{3}+1$ points 
and there exists a $PK$ metric on $\mathbb CP^2$
with conical angles $2\pi\frac{N-3}{N}$ at the lines $L_j$.  
 
\end{corollary}

Note that for a generic arrangement the total multiplicity is $N(N-1)$
while for the most degenerate arrangement it is $N$.

{\bf Proof.} We have the following equality:

$$N^2=(\sum_j L_j)^2=N+\sum_{j\ne k}L_j\cdot L_k=N+\sum_i \mu_i(\mu_i-1)$$

\begin{equation} \label{sumofmu}
\sum_i\mu_i^2=N^2-N+\sum_i\mu_i
\end{equation}
Associate to each $L_j$ weight $\beta_j=\frac{N-3}{N}$, then the 
arrangement satisfies the conditions of Theorem \ref{stabth}.
Since all weights are equal we can treat double points of the 
arrangements as points of type $(1,1)$ and we get the following 
equalities:

$$1-\alpha_i=\frac{3\mu_i}{2N},\;\;\, 1-\beta_j=\frac{3}{N},\;\;\; 
\sum_jB_{jj}=\sum_i\mu_i-N$$

Applying the Bogomolov-Gieseker inequality to the corresponding
stable parabolic bundle $E_*$ and using (\ref{sumofmu}) we get:

$$0\ge{\rm par\textrm{-}ch}_2(E_*)=\frac{\sum_i 9\mu_i^2}{4N^2}-
\frac{9\sum_i \mu_i-9N}{2N^2}-\frac{3}{2}=$$

$$\frac{9\sum_i \mu_i^2-18\sum_i \mu_i+18N-6N^2}{4N^2}=
\frac{3N^2+9N-9\sum_i\mu_i}{4N^2}$$

\hfill$\square$

This inequality was proven previously in [La] using different methods.
We finish with the following remark.

\begin{lemma} For line arrangements that satisfy 
condition of Theorem \ref{th:existence} the system of equations 
(\ref{bgderive}) from Theorem \ref{Thm:chernclasses} is
a corollary of  the Bogomolov-Gieseker inequality. 
\end{lemma}

{\bf Proof.} Let $(L_j,\beta_j)$ be an arrangement satisfying 
conditions of \ref{th:existence}. For a small deformation 
$\beta_j'$ of $\beta_j$ the deformed parabolic bundle $E_*'$
(defined in the same way as $E_*$ but using the weights  
$\beta_j'$ instead $\beta_j$) 
is stable and satisfies the Bogomolov-Gieseker inequality.
The inequality is quadratic in $\beta_j$ and it attains its maxima (zero)
for the initial data. So its derivatives in all directions vanish. This
produces a system of linear equations on $\beta_j$. We will show that 
this system is equivalent to the system of equations (\ref{bgderive})
from Theorem \ref{Thm:chernclasses}. Consider the
following derivative:
$$\frac{\partial}{\partial \beta_j}{\rm par\textrm{-}ch}_{2}(E_*)-
\frac{\partial}{\partial \beta_j}\frac{1}{2}{\rm par\textrm{-}ch}_1^2(E_*)$$ 
If we restrict this expression to the plane
$\sum_j(\beta_j-1)=-3$, the second term vanishes. The first term is 
given by Equation (\ref{parch2cp2})
$$\frac{\partial}{\partial \beta_j}({\rm par\textrm{-}ch}_2(E_*))=
1+(1-\beta_j)B_{jj}-\sum_id_{ij}(1-\alpha_i)=$$

using $\sum_j(\beta_j-1)=-3$ and $1-\alpha_i=\frac{1}{2}\sum_j{d_{ij}(1-\beta_j)}$
$$1+(1-\beta_j)B_{jj}-\frac{1}{2}(3-\sum_{k\ne j}B_{jk}(1-\beta_j)+(B_{jj}+1)(1-\beta_j)))=$$
$$\frac{1}{2}(1-\sum_kB_{jk}(1-\beta_j))$$
\hfill$\square$

\section{Further results, questions, and directions}

In a subsequent paper we will use Theorem \ref{th:existence}
to construct several infinite families of aspherical complex
surfaces. Some of these families of surfaces
admit a metric of type $CAT(0)$. In particular some 
smooth compact quotients of the complex ball 
admit a $PK$ metric of type $CAT(0)$, 
this answers a question of Gromov and Davis.

For every $PK$ surface its $PK$ metric induces a positive $(1,1)$
current on it, so it should  not be difficult to prove that every 
$PK$ surface is a K\"ahler surface
(in principle, one should be able to smoothen a bit the $PK$
metric to get a smooth K\"ahler metric on the surface). 
It should be possible to show that
non-algebraic $K3$ surfaces 
don't admit a $PK$ metric, but we don't know any 
obstruction for the existence of $PK$ metrics on algebraic surfaces. 
At the same time the set of examples of $PK$ surfaces 
that we have is rather limited.

We hope that  Theorem \ref{pkcomplex} extends 
to higher dimensions. Namely, that for a polyhedral 
K\"ahler manifold of any dimension the complex structure 
on the complement of the metric singularities extends to a complex 
analytic structure on the whole manifold. In particular
the metric singularities should not have odd (real) codimension
and all singularities of even codimension should have holomorphic 
directions. Note that in higher dimensions 
we can obtain  complex manifolds with singularities  
even if we start with a topological polyhedral 
K\"ahler manifold (i.e. the link of every point is a 
topological sphere). Indeed, by Brieskorn the link of the hypersurfaces
$z_1^2+z_2^2+z_3^2+z_4^3+z_5^{6k-1}=0$, $1 \le k\le 28$ in $\mathbb C^5$
is $S^7$ with one of $28$ smooth structures. 
At the same time these hypersurfaces have a $PK$ metric, 
induced by an obvious degree $24(6k-1)$ 
ramified cover of the hyperplane $\sum_i z_i=0$.

The notion of polyhedral K\"ahler manifolds can be 
generalised in several directions. A 
{\it polyhedral affine structure} on a manifold is a choice of 
a simplicial decomposition and an affine structure on the complement of 
codimension $2$ faces that restricts to the standard affine structure on
the interior of every face of the top dimension. 
We say that a manifold $M^{2n}$ is 
{\it polyhedral complex affine} if the holonomy is contained 
in $GL(n,\mathbb C)$ and singular faces of codimension $2$ at which
the holonomy is trivial have 
holomorphic directions ({\it cf} Definition \ref{gooddefinition}). 
For complex dimension $2$
we expect to get a theory similar to the one developed in this 
article. It should be possible to classify the singularities 
of complex codimension $2$ but the list will be longer.  

If the holonomy of a polyhedral affine manifold is contained in 
the symplectic group $SP(2n)$ we call the manifold 
{\it polyhedral nearly symplectic}. It is not hard to see 
that every symplectic manifold admits a polyhedral nearly  
symplectic structure. But the converse 
should be wrong already for $4$-manifolds, so we adjust the 
definition.

\begin{definition} 
Let $M^4$ be a polyhedral nearly symplectic manifold and let $M_1^4(\varepsilon)$
be a neighborhood of the union of all edges. The 
$PL$ symplectic structure on 
$M^4\setminus M_1^4(\varepsilon)$ can be smoothen 
along the faces of codimension $2$ to a genuine symplectic form $w$. 
Let $c_1$ be the first Chern class of an almost 
complex structure tamed by $w$ on $M_1^4(\varepsilon)$. We call $M^4$
{\it polyhedral symplectic} if for every surface $S$ contained in 
$M_1^4(2\varepsilon)\setminus M_1^4(\varepsilon)$ we have $c_1\cdot S=0$.
\end{definition}
 
It is not hard to prove that every symplectic $4$-manifold
admits a polyhedral symplectic structure. More importantly, 
we conjecture that every polyhedral symplectic $4$-manifold 
admits a symplectic smoothing.






Another interesting direction to generalise $PK$ manifolds is to 
consider complex manifolds with K\"ahler metric of constant 
holomorphic curvature and conical singularities at 
holomorphic geodesic divisors. In the case of negative 
curvature these manifolds will be 
generalisations of complex hyperbolic orbifolds, i.e. 
quotients of  the complex ball $B^n$ by a lattice of 
$SU(n,1)$. These manifolds are presumably are the 
same as Thurston $(\mathbb CH^n, SU(n,1))$ - cone manifolds [Th].
In the case of surfaces it is sufficient 
to ask that singularities of the metric 
form a complex curve, at points that are not 
multiple there is a local isometric action of the 
group $SU(1,1)\times S^1$, and at the multiple points 
there is an holomorphic isometric action of $\mathbb R^1$. 
It should be possible to generalise Theorem  \ref{th:existence}
to this setting using parabolic Kobayashi-Hitchin correspondence 
for parabolic Higgs bundles [M2] using an idea from [S].
It would be interesting to reprove (or even generalise)
results of [CHL] using this approach.

Finally we hope to prove the following conjecture.

\begin{conjecture}
 For every arrangement satisfying conditions 
of Theorem \ref{th:existence} its complement is of the type $K(\pi,1)$.
\end{conjecture}

The converse to this conjecture is wrong because 
Proposition \ref{doublecrit} permits us to check that 
some  simplicial arrangements are not $PK$. At the same time 
by a theorem of Deligne all (complexified) simplicial arrangement 
have a complement of the type $K(\pi,1)$.


\begin{thebibliography}{99999999999999}

\bibitem[Bo]{Bo} Borel et al. 
Algebraic $D$-modules. Perspectives in Mathematics, 2.
(1987) 151--172.

\bibitem[BK]{BK} T. F. Banchoff, W. K{\"u}hnel.
 The $9$-vertex complex projective plane.
 {\it Math. Intelligencer.}  {\bf 5}(1983), 11--22.

\bibitem[Ch]{Ch}J. Cheeger.
A vanishing theorem for piecewise constant curvature spaces.
{\it Curvature and topology of Riemannian manifolds (Katata, 1985), 33--40,
Lecture Notes in Math., 1201,
Springer, Berlin.} (1986).

\bibitem[CHL]{CHL}  W. Couwenberg, G. Heckman, E. Looijenga.
 Geometric structures on the complement of a projective arrangement.
{\it Publ. Math. IHES.} {\bf 101} (2005), 69--161.

\bibitem[F]{F} R. Friedman. Algebraic surfaces and holomorphic
vector bundles. Springer (1998).

\bibitem[G1]{G1} B. Gr\"unbaum. Arrangements of hyperplanes.
Proceedings of the second Louisisana conference on combinatorics 
graph theory and computings. 41--106 (1971).

\bibitem[G2]{G2} B. Gr\"unbaum. Arrangements and spreads. Regional
Conf. Series in Mathematics {\bf 10}. Amer. Math Soc., (1972).   

\bibitem[Hir]{Hir}F. Hirzebruch. Algebraic surfaces with extreme Chern numbers.
{\it Russian Math. Surveys} {\bf 40} (1985), 135--145.

\bibitem[IS]{IS} Jaya NN Iyer, C. Simpson. The Chern character of
a parabolic bundle, and a parabolic Reznikov theorem in the case of
finite order at infinity. arXiv: math.AG/06121444


\bibitem[KTM]{KTM} J. Kaneko, S. Tokungaga, M. Yoshida.
Complex crystallographic
groups 2. {\it J. Math. Soc. Japan.}  {\bf 34}(1982), 595--605.


\bibitem[La]{La} A. Langer. Logarithmic orbifold Euler numbers of 
surfaces with applications. {\it Proc. London Math. Soc.} 
(3) 86 (2003), no. 2, 358--396. 

\bibitem[Li]{Li} Li Jiayu. Hermitian-Einstein metrics and Chern number 
inequalities on parabolic stable bundles over K\"ahler manifolds.
{\it Communications in analysis and geometry.} {8}(2000), 445--475   



\bibitem[M1]{M1} T. Mochizuki. Asymptotic behaviour of tame harmonic 
bundles and an application to pure twistor $D$-modules, part 1. 
{\it Mem. Amer. Math. Soc.} 185 (2007), no. 869, xii+324 pp

\bibitem[M2]{M2} T. Mochizuki. Kobayashi-Hitchin correspondence for tame
harmonic bundles and an application. {\it Ast\'erisque} No. 309 (2006)


\bibitem[Oh]{Oh} M. Ohtsuki. A residue formula for Chern classes
associated with logarithmic connections.{\it Tokyo J. Math.}
{\bf  5}(1982) 13--21.

\bibitem[Or]{Or} S. Orshanskiy. A PL-manifold of 
nonnegative curvature homeomorphic to $S^2 \times S^2$ 
is a direct metric product.   arXiv:0807.1922 (2008).

\bibitem[S]{S}C. Simpson. Constructing variations of Hodge 
structure using Yang-Mills theory
and application to uniformization. {\it J. Amer. Math. Soc.} 1 (1988), 867--918.

\bibitem[Th]{Th} W. P. Thurston. Shapes of polyhedra and triangulations
of the sphere.
{\it Geometry and Topology Monographs} {\bf 1}, (1998) 511--549.

\bibitem[Tr]{Tr} M. Troyanov. Les surfaces euclidiennes {\`a}
singularit{\'e}s coniques.
{\it L'Enseign. Math.} {\bf 32}(1986), 79--94.

\end{thebibliography}
\end{document}